\documentclass[twoside]{irmaems}
\usepackage{amssymb} 
\usepackage{amsmath} 
\usepackage{latexsym}

\usepackage[mathcal]{eucal}
\usepackage{graphicx}
\usepackage{pinlabel}
\usepackage{amscd}

\renewcommand\theequation{\arabic{section}.\arabic{equation}}

\theoremstyle{definition} 

 \newtheorem{definition}{Definition}[section]

\theoremstyle{plain}      

 \newtheorem{prop}[definition]{Proposition}
 \newtheorem{theorem}[definition]{Theorem}
 \newtheorem{corollary}[definition]{Corollary}
 \newtheorem{lemma}[definition]{Lemma}

\begin{document}

\title{A survey of quantum Teichm\"uller space and Kashaev algebra}

\author{Ren Guo}

\address{
School of Mathematics\\
University of Minnesota\\
Minneapolis, MN 55455\\
email:\, \tt{guoxx170@math.umn.edu}}

\maketitle

\begin{abstract} 

In this chapter, we survey the algebraic aspects of quantum Teichm\"uller space, generalized Kashaev algebra and a natural relationship between the two algebras.

\end{abstract}

\begin{classification}
57R56; 57M50, 20G42.
\end{classification}

\begin{keywords}
shear coordinate, quantum Teichm\"uller space, Kashaev coordinates, generalized Kashaev algebra.
\end{keywords}

\tableofcontents

\section{Introduction}

A quantization of the Teichm\"uller space $\mathcal T(S)$\index{quantization of Teichm\"uller spaces} of a punctured surface $S$ was developed by L. O. Chekhov and V. V. Fock \cite{FC,CF,Fo} and, independently, by R. Kashaev \cite{Kas1, Kas2, Kas3, Kas4}. This is a deformation of the $\mathrm C^*$--algebra of functions on Teichm\"uller space $\mathcal T(S)$. The quantization was expressed in terms of self-adjoint operators on Hilbert spaces and the quantum dilogarithm function. Although these two approaches of quantization use the same ingredients, the relationship between them is still mysterious. Chekhov and Fock worked with shear coordinates of Teichm\"uller space while Kashaev worked with a new coordinate.

The pure algebraic foundation of Chekhov-Fock's quantization was established in the work of F. Bonahon, H. Bai and X. Liu \cite{BBL, BonLiu, Liu1}. The algebraic aspect of Kashaev's quantization is investigated and generalized in \cite{GL}. And a natural relationship between quantum Teichm\"uller space and generalized Kashaev algebra is established in \cite{GL}. In this chapter we make a survey of the ideas and results mentioned.

Recently, I. B. Frenkel and H. K. Kim \cite{FK} derived the quantum Teichm\"uller space from tensor products of a single canonical representation of the modular double of the quantum plane and showed that the quantum universal Teichm\"uller space is realized in the infinite tensor power of the canonical representation naturally indexed by rational numbers including the infinity.

\noindent \textbf{Acknowledgment.} We would like to thank Athanase Papadopoulos for inviting us to contribute a chapter to the series of Handbook of Teichm\"uller space, for careful reading of the manuscript and for suggestions on improving the exposition of the chapter. We also thank
the referee for helpful suggestions to improve the chapter.

\section{The quantum Teichm\"uller space: foundation}

In this section we review the finite-dimensional Chekhov-Fock's quantization of the Teichm\"uller space following \cite{Liu1} closely. 

\subsection{Ideal triangulations}

Let $S$ be an oriented surface with genus $g$ and with $p\geq 1$ punctures, obtained by removing $p$ points $\{v_1,\ldots,v_p\}$ from a closed oriented surface $\bar S$ of genus $g$. An \emph{ideal triangulation}\index{ideal triangulation} of $S$ is a triangulation of the closed surface $\bar{S}$ whose vertex set is exactly $\{v_1,\ldots,v_p\}$. If the Euler characteristic of $S$ is negative, i.e., $m:=2g-2+p>0,$ $S$ has an ideal triangulation. Any ideal triangulation of $S$ has $2m$ ideal triangles and $3m$ edges. The edges of an ideal triangulation $\lambda$ of $S$ are numerated as $\{\lambda_1,...,\lambda_{3m}\}$.

Let $\Lambda(S)$ denote the set of isotopy classes of ideal triangulations of $S$. The set $\Lambda(S)$ admits a natural action of the symmetric group on the set $\{ 1, 2, ..., 3m \}$, $\mathfrak{S}_{3m}$, acting by permuting the indices of the edges of $\lambda$. Namely $\lambda'
= \alpha(\lambda)$ for $\alpha \in\mathfrak{S}_{3m}$ if $\lambda_i=\lambda'_{\alpha(i)}$.

\begin{figure}[ht!]
\labellist\small\hair 2pt
\pinlabel $\lambda_i$ at 76 57
\pinlabel $\lambda'_i$ at 381 78
\pinlabel $\longrightarrow$ at 216 60
\pinlabel $\Delta_i$ at 216 73

\endlabellist
\centering
\includegraphics[scale=0.5]{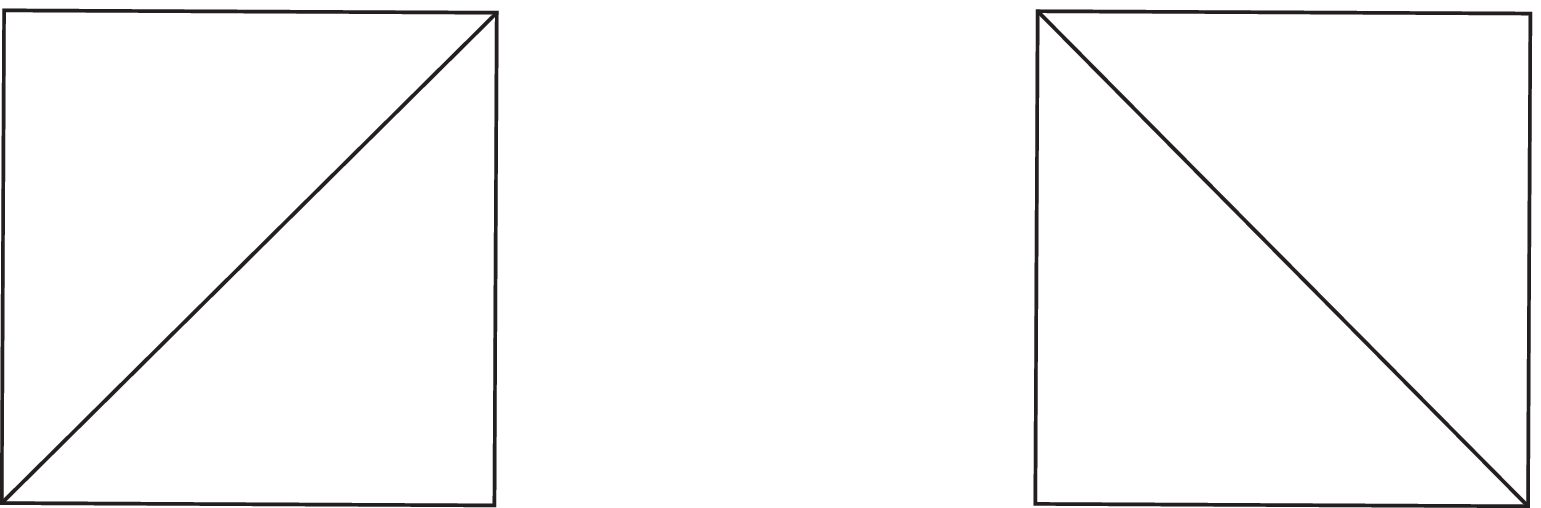}
\caption{}
\label{fig:diagonal-exchange}
\end{figure}

Another important transformation of $\Lambda(S)$ is provided by the \emph {$i$--th diagonal exchange map}\index{diagonal exchange} $\Delta_i: \Lambda(S)\rightarrow \Lambda(S)$ defined as follows. Suppose that the $i$--th edge $\lambda_i$ of an ideal triangulation $\lambda\in \Lambda(S)$ is adjacent to two triangles. Then $\Delta_i(\lambda)$ is obtained from $\lambda$ by replacing the edge $\lambda_i$ by the other diagonal $\lambda_i'$ of the square formed by the two triangles, as illustrated in Figure \ref{fig:diagonal-exchange}.

\begin{figure}[ht!]
\labellist\small\hair 2pt
\pinlabel $\lambda_i$ at 46 83
\pinlabel $\lambda_j$ at 120 83

\endlabellist
\centering
\includegraphics[scale=0.5]{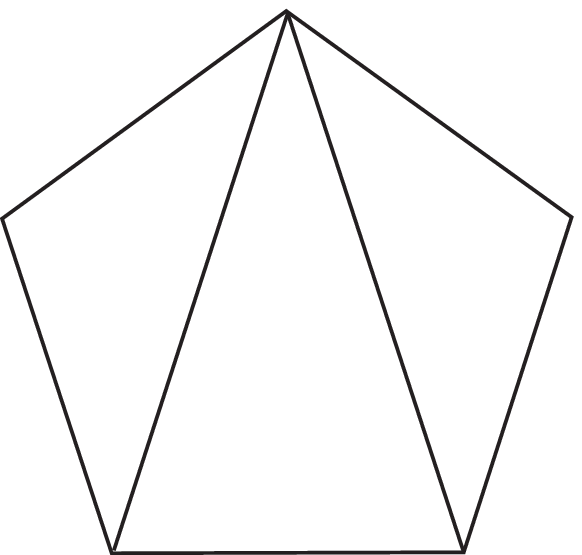}
\caption{}
\label{fig:pentagon-ideal}
\end{figure}

\begin{lemma}\label{thm:graph-ideal} 
The reindexings and diagonal exchanges satisfy the following relations:

\begin{enumerate}
\item $(\alpha\beta)(\lambda) = \alpha(\beta(\lambda))$ for every $\alpha$, $\beta\in\mathfrak{S}_{3m}$;

\item $(\Delta_i)^2=\mathrm{Id}$;

\item $\alpha  \circ \Delta_i = \Delta_{\alpha(i)}  \circ \alpha$ for every $\alpha\in\mathfrak{S}_{3m}$;

\item If $\lambda_i$ and $\lambda_j$ do not belong to the same triangle of $\lambda\in \Lambda(S)$, then $\Delta_i \circ \Delta_j(\lambda)= \Delta_j
\circ \Delta_i(\lambda)$;

\item If three triangles of an ideal triangulation $\lambda\in \Lambda(S)$ form a pentagon with
diagonals $\lambda_i$, $\lambda_j$ as in Figure \ref{fig:pentagon-ideal}, then
\begin{equation*}
\Delta_i\circ \Delta_j\circ \Delta_i\circ \Delta_j\circ \Delta_i(\lambda) = \alpha_{i\leftrightarrow j}(\lambda),
\end{equation*}
where $\alpha_{i\leftrightarrow j} \in\mathfrak{S}_{3m}$ denotes the
transposition exchanging $i$ and $j$.

\end{enumerate}

\end{lemma}

To construct the quantum Teichm\"uller space, we need the following two results of R. C. Penner \cite{Pen} (see also J. L. Harer \cite{Har}).

\begin{theorem}\label{thm:Penner1-ideal} 
Given  two ideal triangulations $\lambda,\lambda'\in \Lambda(S)$, there exists a finite sequence of ideal triangulations $\lambda=\lambda_{(0)}$, $\lambda_{(1)}$, \dots, $\lambda_{(n)}=\lambda'$ such that each $\lambda_{(k+1)}$ is obtained from $\lambda_{(k)}$ by a diagonal exchange or by a
reindexing of its edges. 
\end{theorem}

\begin{theorem}\label{thm:Penner2-ideal} 
Given  two ideal triangulations
$\lambda,\lambda'\in \Lambda(S)$ and given two sequences
$\lambda=\lambda_{(0)}$, $\lambda_{(1)}$, \ldots, $\lambda_{(n)}
=\lambda'$ and $\lambda=\overline{\lambda}_{(0)}$, $\overline{\lambda}_{(1)}$, \ldots,
$\overline{\lambda}_{(\overline{n})} =\lambda'$ of diagonal exchanges and reindexings
connecting them as in Theorem~\ref{thm:Penner1-ideal}, these two
sequences can be related to each other by successive applications
of the following moves and of their inverses. These moves correspond to the relations in Lemma \ref{thm:graph-ideal}.
\begin{enumerate}

\item Replace \dots, $\lambda_{(k)}$, $\beta(\lambda_{(k)})$, $\alpha( \beta(\lambda_{(k)}))$, \dots\ 

\hspace{22pt}by \dots, $\lambda_{(k)}$, $(\alpha \beta)(\lambda_{(k)})$, \dots\ where $\alpha$, $\beta\in \mathfrak{S}_{3m}$.

\item Replace \dots, $ \lambda_{(k)}$, $\Delta_i(\lambda_{(k)}) $, $\lambda_{(k)}$, \dots\ 

\hspace{22pt}by \dots, $\lambda_{(k)}$, \dots\ .

\item Replace \dots, $\lambda_{(k)}$, $\Delta_i(\lambda_{(k)})$, $  \alpha \circ \Delta_i(\lambda_{(k)})$, \dots\ 

\hspace{22pt}by \dots , $\lambda_{(k)}$, $\alpha(\lambda_{(k)})$, $ \Delta_{\alpha(i)}\circ \alpha(\lambda_{(k)})$, \dots\ where
$\alpha\in \mathfrak{S}_{3m}$.

\item Replace  \dots, $\lambda_{(k)}$, $\Delta_i(\lambda_{(k)})$, $\Delta_j\circ \Delta_i(\lambda_{(k)})$, \dots\ 

\hspace{22pt}by \dots, $ \lambda_{(k)}$, $\Delta_j(\lambda_{(k)})$, $\Delta_i\circ \Delta_j(\lambda_{(k)})$, \dots\  where
$\lambda_i,\lambda_j$ are two edges which do not belong to a same
triangle of $\lambda_{(k)}$.

\item Replace \dots, $\lambda_{(k)}$, $\Delta_i(\lambda_{(k)})$,
$\Delta_j\circ \Delta_i(\lambda_{(k)})$, $\Delta_i\circ
\Delta_j\circ \Delta_i(\lambda_{(k)})$, $\Delta_j\circ
\Delta_i\circ \Delta_j\circ \Delta_i(\lambda_{(k)})$, $\Delta_i\circ\Delta_j\circ
\Delta_i\circ \Delta_j\circ \Delta_i(\lambda_{(k)})$, \dots\ 

\hspace{22pt}by \dots, $\lambda_{(k)}$,
$\alpha_{i\leftrightarrow j} (\lambda_{(k)})$, \dots\ where $\lambda_i,\lambda_j$ are two diagonals of a pentagon
of $\lambda_{(k)}$ as in Figure \ref{fig:pentagon-ideal}.

\end{enumerate}
\end{theorem}

\subsection{Shear coordinates for the Teichm\"uller space}

If the Euler characteristic of $S$ is negative, i.e., $m:=2g-2+p>0,$ $S$ admits complete hyperbolic metrics. The \emph{Teichm\"uller space}\index{Teichm\"uller space} $\mathcal T(S)$ of $S$ consists of all isotopy classes of complete hyperbolic metrics on $S$. 
W. Thurston \cite{Th} associated to each ideal triangulation a global coordinate system which is called \emph{shear coordinate}\index{shear coordinate} for the Teichm\"uller space $\mathcal T(S)$ (see also \cite{Bon, Fo}).

An end of a surface $S$ with a complete hyperbolic metric $d \in \mathcal{T}(S)$ can be of two types: a \emph{cusp}\index{cusp} with finite area bounded on
one side by a horocycle; and a \emph{funnel}\index{funnel} with infinite area bounded
on one side by a simple closed geodesic. The \emph{convex core}\index{convex core} $\mathrm{Conv}(S,d)$ of $(S, d)$ is the smallest non-empty closed convex subset
of $(S,d)$, and is bounded in $S$ by a family of disjoint simple closed geodesics. Each cusp end of $(S,d)$ is also a cusp end of $\mathrm{Conv}(S,d)$, while each funnel end of $S$ faces a boundary component of $\mathrm{Conv}(S,d)$.

The \emph{enhanced Teichm\"uller space}\index{Teichm\"uller space!enhanced} $\widetilde{\mathcal
T}(S)$ consists of all isotopy classes of complete hyperbolic metrics $d\in
\mathcal T(S)$ enhanced with an orientation of each boundary component of $\mathrm
{Conv}(S,d)$. 

Under an enhanced hyperbolic metric $d \in \widetilde{\mathcal
T}(S)$, each edge
$\lambda_i$ of an ideal triangulation $\lambda$ is realized by a unique
$d$--geodesic $g_i$ such that each end of $g_i$, either converges
towards a cusp end of $S$,  or spirals around a boundary component of
$\mathrm{Conv}(S,d)$ in the orientation specified by $d\in \widetilde{\mathcal T}(S)$. 

\begin{figure}[ht!]
\labellist\small\hair 2pt
\pinlabel $\mathbb H^2$ at 351 166
\pinlabel $\widetilde g_i$ at 251 137
\pinlabel $\widetilde T_i^1$ at 167 113
\pinlabel $\widetilde T_i^2$ at 233 82
\pinlabel $\widetilde Q$ at 99 150
\pinlabel $z_-$ at 179 21
\pinlabel $z_+$ at 358 21
\pinlabel $z_{\mathrm r}$ at 239 21
\pinlabel $z_{\mathrm l}$ at 81 21

\endlabellist
\centering
\includegraphics[scale=0.55]{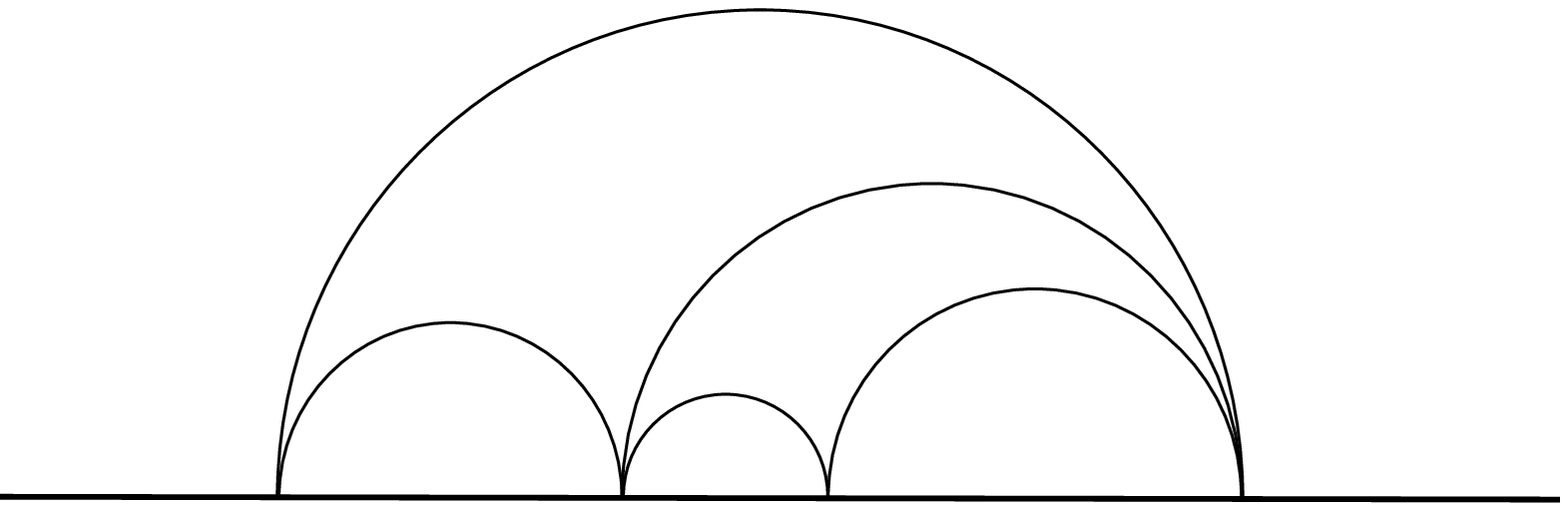}
\caption{}
\label{fig:shear}
\end{figure}

The enhanced hyperbolic metric $d \in \widetilde{\mathcal T}(S)$
associates to the edge $\lambda_i$ of $\lambda$ a positive
number $x_i$ defined as follows. The geodesic $g_i$ separates two
triangle components $T_i^1$ and $T_i^2$ of $\mathrm{Conv}(S,d)-\{g_i\}$.
The hyperbolic plane $\mathbb H^2$ is the universal covering of $S$ endowed with the metric $d$. Lift $g_i$, $T_i^1$ and
$T_i^2$ to a geodesic $\widetilde g_i$ and two triangles
$\widetilde T_i^1$ and $\widetilde T_i^2$ in $\mathbb H^2$ so
that the union $\widetilde g_i \cup \widetilde T_i^1\cup
\widetilde T_i^2$ forms a square $\widetilde Q$ in $\mathbb H^2$. See Figure \ref{fig:shear}.
In the upper half-space model for $\mathbb H^2$, let
$z_-$, $z_+$, $z_{\mathrm r}$, $z_{\mathrm l}$ be the vertices of
$\widetilde Q$ in such a way that $\widetilde g_i$ goes from $z_-$
to $z_+$ and, for this orientation of $\widetilde g_i$,
$z_{\mathrm r}$, $z_{\mathrm l}$ are respectively to the right and
to the left of $\widetilde g_i$ for the orientation of $\widetilde
Q$ given by the orientation of $S$. Then,
$$x_i:=- \,\textrm{cross-ratio}\,
(z_{\mathrm r},z_{\mathrm l},z_-,z_+)=-\, \frac{(z_{\mathrm
r}-z_-)(z_{\mathrm l}-z_+)}{(z_{\mathrm r}-z_+) (z_{\mathrm
l}-z_-)}.$$ 

The real numbers $\{x_i\}$ are the \emph{exponential shear
coordinates}\index{shear
coordinates!exponential} of the enhanced hyperbolic metric $d \in
\widetilde{\mathcal T}(S)$. The shear coordinates are $\ln x_i$.

It turns out that $\{x_i\}$ defines a homeomorphism $\phi_\lambda
:\widetilde{\mathcal T}(S) \rightarrow \mathbb R_+^{3m}$.

Therefore the exponential shear coordinates associates a parametrization
$\phi_\lambda :\widetilde{\mathcal T}(S) \rightarrow \mathbb
R_+^{3m}$ to each ideal triangulation $\lambda\in \Lambda(S)$
(endowed with an indexing of its edges). We now investigate the
coordinate changes $\phi_{\lambda'} \circ \phi_{\lambda}^{-1}$
associated to two ideal triangulations.

If $\lambda'= \alpha(\lambda)$ is obtained by reindexing the edges
of $\lambda$ by $\alpha \in \mathfrak S_{3m}$, then
$\phi_{\lambda'} \circ \phi_{\lambda}^{-1}$ is the
permutation of the coordinates by $\alpha$. For a diagonal
exchange, we have the following result.
 
\begin{figure}[ht!]
\labellist\small\hair 2pt
\pinlabel $\lambda_i$ at 76 57
\pinlabel $\lambda_j$ at 66 159
\pinlabel $\lambda_m$ at -13 70
\pinlabel $\lambda_l$ at 66 -15
\pinlabel $\lambda_k$ at 155 70
\pinlabel $\lambda'_i$ at 381 78
\pinlabel $\lambda'_j$ at 365 159
\pinlabel $\lambda'_m$ at 285 70
\pinlabel $\lambda'_l$ at 365 -15
\pinlabel $\lambda'_k$ at 453 70
\pinlabel $\longrightarrow$ at 216 60
\pinlabel $\Delta_i$ at 216 73

\endlabellist
\centering
\includegraphics[scale=0.5]{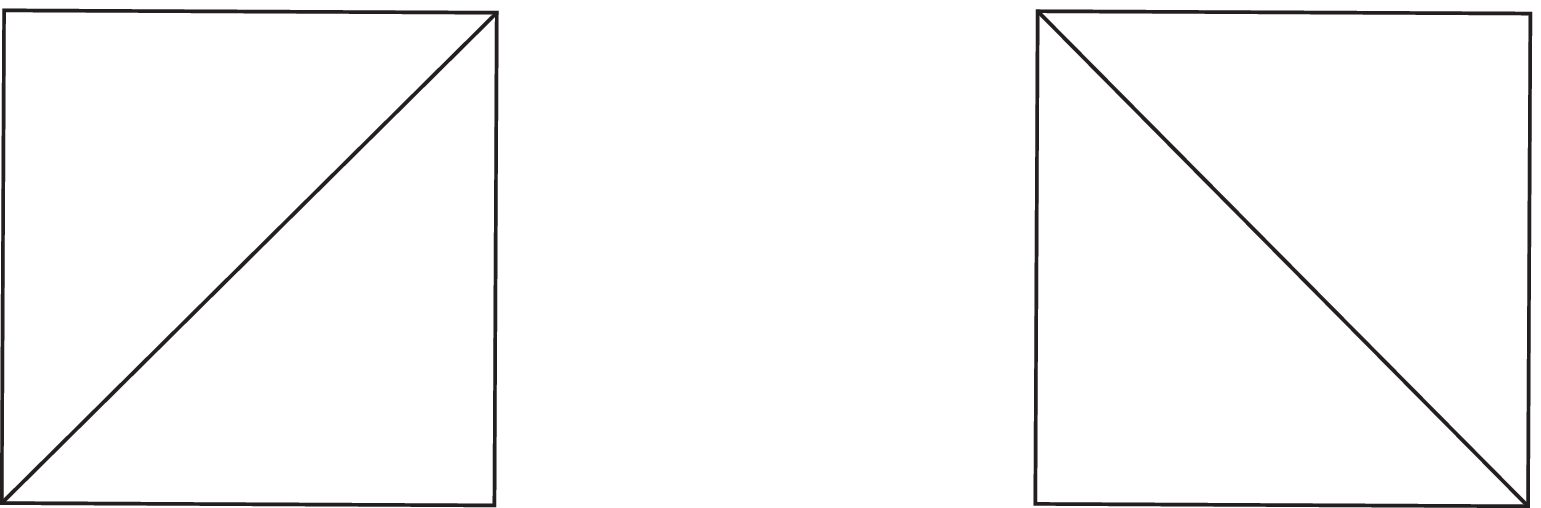}
\caption{}
\label{fig:eight-ideal}
\end{figure} 
 
\begin{prop}[Liu \cite{Liu1}]
\label{prop:DiagExch} Suppose that the ideal triangulations
$\lambda$, $\lambda'$ are obtained from each other
by a diagonal exchange, namely that $\lambda' =
\Delta_i(\lambda)$. Label the edges of $\lambda$ involved in this
diagonal exchange as $\lambda_i$, $\lambda_j$, $\lambda_k$,
$\lambda_l$, $\lambda_m$ as in Figure~\ref {fig:eight-ideal}. If
$(x_1, x_2, \dots, x_{3m})$ and $(x_1', x_2', \dots,
x_{3m}')$ are the exponential shear coordinates
associated to $\lambda$ and $\lambda'$ of the same enhanced hyperbolic metric,
then $x_h'=x_h$ for every $h \not\in
\{i,j,k,l,m\}$, $x_i' = x_i^{-1}$ and:
\begin{description}
\item[Case 1] if the edges $\lambda_j$, $\lambda_k$, $\lambda_l$,
$\lambda_m$  are distinct, then
$$ x'_j = (1+x_i)x_j\quad x'_k  =
(1+x_i^{-1})^{-1}x_k \quad x'_l = (1+x_i)x_l \quad x'_m =
(1+x_i^{-1})^{-1}x_m;$$

\item[Case 2] if $\lambda_j$ is identified with $\lambda_k$, and
$\lambda_l$ is distinct from $\lambda_m$, then
$$
    x'_j  =  x_ix_j \quad x'_l  =
(1+x_i)x_l \quad x'_m  =  (1+x_i^{-1})^{-1}x_m;
$$

\item[Case 3] (the inverse of Case 2) if $\lambda_j$ is identified
with $\lambda_m$, and $\lambda_k$ is distinct from $\lambda_l$,
then
$$
    x'_j  =  x_ix_j
\quad x'_k  =  (1+x_i^{-1})^{-1}x_k \quad x'_l  = (1+x_i)x_l;
$$

\item[Case 4] if $\lambda_j$ is identified with $\lambda_l$, and
$\lambda_k$ is distinct from $\lambda_m$, then
$$
    x'_j  =  (1+x_i)^2x_j
\quad x'_k  =  (1+x_i^{-1})^{-1}x_k \quad x'_m  =
(1+x_i^{-1})^{-1}x_m
$$

\item[Case 5] (the inverse of Case 4) if $\lambda_k$ is identified
with $\lambda_m$, and $\lambda_j$ is distinct from $\lambda_l$,
then
$$ x'_j  =  (1+x_i)x_j \quad x'_k  =
(1+x_i^{-1})^{-2}x_k \quad x'_l = (1+x_i)x_l;
$$

\item[Case 6] if $\lambda_j$ is identified with $\lambda_k$, and
$\lambda_l$ is identified with $\lambda_m$ (in which case $S$ is a
$3$--times punctured sphere), then
$$
    x'_j  =  x_ix_j \quad x'_l  =
x_ix_l;
$$

\item[Case 7] (the inverse of Case 6) if $\lambda_j$ is identified
with $\lambda_m$, and $\lambda_k$ is identified with $\lambda_l$
(in which case $S$ is a $3$--times punctured sphere), then
$$
    x'_j  =  x_ix_j
\quad x'_k  =  x_ix_k;
$$

\item[Case 8] if $\lambda_j$ is identified with $\lambda_l$, and
$\lambda_k$ is identified with $\lambda_m$ (in which case $S$ is a
once punctured torus), then
$$ x'_j  =  (1+x_i)^2x_j \quad x'_k  =
(1+x_i^{-1})^{-2}x_k.
$$

\end{description}
\end{prop}

\subsection{The Chekhov-Fock algebra}

Fix an ideal triangulation $\lambda\in \Lambda(S)$. The complement
$S-\lambda$ has $6m$ spikes converging towards the punctures, and
each spike is delimited by one $\lambda_i$ on one side and one
$\lambda_j$ on the other side, with possibly $i=j$. For $i$, $j\in
\{1, \dots, 3m\}$, let $a^{\lambda}_{ij}$ denote the number of
spikes of $S-\lambda$ which are delimited on the left by
$\lambda_i$ and on the right by $\lambda_j$, and set
\begin{equation*}
\sigma^{\lambda}_{ij}= a^{\lambda}_{ij}-a^{\lambda}_{ji}.
\end{equation*} Note that
$\sigma_{ij}^\lambda \in \{-2, -1, 0, 1, 2\}$, and that $\sigma_{ji}^\lambda= -\sigma_{ij}^\lambda$.

Let $q$ be an arbitrary complex number. The \emph{Chekhov-Fock algebra}\index{Chekhov-Fock algebra} associated to the ideal
triangulation $\lambda$ is the algebra $\mathcal{T}^q_{\lambda}$
defined by generators $X_1$, $X_1^{-1}$, $X_2$, $X_2^{-1}$, \dots,
$X_{3m}$, $X_{3m}^{-1}$, with each pair $X_i^{\pm 1}$ associated to an
edge $\lambda_i$ of $\lambda$, and by the relations
\begin{equation*}
X_iX_j=q^{2\sigma^{\lambda}_{ij}}X_jX_i.
\end{equation*}

This algebra has a well-defined fraction division algebra $\widehat{\mathcal T}_\lambda^q$ which
consists of all formal fractions $PQ^{-1}$ with $P$, $Q\in
\mathcal{T}^q_{\lambda}$ and $Q\neq 0$, and two such fractions
$P_1Q_1^{-1}$ and $P_2Q_2^{-1}$ are identified if there exists
$S_1$, $S_2\in \mathcal{T}^q_{\lambda} - \{0\}$ such that $P_1S_1
= P_2S_2$ and $Q_1S_1 = Q_2S_2$.

The algebras ${\mathcal T}_\lambda^q$ and
$\widehat{\mathcal{T}}^q_{\lambda}$ strongly depend on the ideal
triangulation $\lambda$. As one moves from one ideal triangulation $\lambda$  to another $\lambda'$, Chekhov and Fock \cite{Fo, FC, CF} (see also \cite{Liu1}) introduce \emph{coordinate change isomorphisms}\index{coordinate change isomorphism} $\Phi_{\lambda\lambda'}^q: \widehat{\mathcal T}_{\lambda'} ^q \rightarrow \widehat{\mathcal T}_\lambda^q$. 
We denote by $X'_1$, $X_2'$, \dots,
$X_n'$ the generators of $\widehat{\mathcal{T}}^q_{\lambda'}$
associated to the edges $\lambda_1'$, $\lambda_2'$, \dots,
$\lambda_n'$ of $\lambda'$, and by $X_1$, $X_2$, \dots, $X_n$ the
generators of $\widehat{\mathcal{T}}^q_{\lambda}$ associated to
the edges $\lambda_1$, $\lambda_2$, \dots, $\lambda_n$ of
$\lambda$.

\begin{definition}\label{def:exchange} 
Suppose that the ideal triangulations
$\lambda$, $\lambda'\in \Lambda(S)$ are obtained from each other
by an edge reindexing, namely that $\lambda'_i =
\lambda_{\alpha(i)}$ for some permutation $\alpha \in \mathfrak
S_{3m}$.  Then we define a map $\widehat{\alpha}$ from the set of the generators of the algebra $\widehat{\mathcal{T}}^q_{\lambda'}$ to $\widehat{\mathcal{T}}^q_{\lambda}$
by $$\widehat{\alpha}(X'_i) = X_{\alpha(i)}, \ \ \mbox{for any}\ \ i=1,...,3m.$$

Suppose that the ideal triangulations
$\lambda$, $\lambda'$ are obtained from each other
by a diagonal exchange, namely that $\lambda' =
\Delta_i(\lambda)$. Label the edges of $\lambda$ involved in this
diagonal exchange as $\lambda_i$, $\lambda_j$, $\lambda_k$,
$\lambda_l$, $\lambda_m$ as in Figure \ref {fig:eight-ideal}. Then
we define a map $\widehat{\Delta}_i$ on the set of the generators of the algebra $\widehat{\mathcal{T}}^q_{\lambda'}$ to $\widehat{\mathcal{T}}^q_{\lambda}$
such that
    $X_h' \mapsto X_h$
for every $h \not\in \{i,j,k,l,m\}$, $X_i' \mapsto X_i^{-1}$ and:
\begin{description}
\item[Case 1] if the edges $\lambda_j$, $\lambda_k$, $\lambda_l$,
$\lambda_m$  are distinct, then
\begin{align*}
X'_j &\mapsto (1+qX_i)X_j \qquad\, X'_k \mapsto
(1+qX_i^{-1})^{-1}X_k \\
    X'_l &\mapsto (1+qX_i)X_l
\qquad
    X'_m
\mapsto (1+qX_i^{-1})^{-1}X_m ;
\end{align*}

\item[Case 2] if $\lambda_j$ is identified with $\lambda_k$, and
$\lambda_l$ is distinct from $\lambda_m$, then
$$ X'_j  \mapsto   X_iX_j \quad X'_l
\mapsto (1+qX_i)X_l \quad X'_m \mapsto (1+qX_i^{-1})^{-1}X_m
$$

\item[Case 3] (the inverse of Case 2) if $\lambda_j$ is identified
with $\lambda_m$, and $\lambda_k$ is distinct from $\lambda_l$,
then
$$X'_j\mapsto   X_iX_j \quad X'_k
\mapsto (1+qX_i^{-1})^{-1}X_k \quad X'_l  \mapsto (1+qX_i)X_l
$$

\item[Case 4] if $\lambda_j$ is identified with $\lambda_l$, and
$\lambda_k$ is distinct from $\lambda_m$, then
\begin{gather*}
X'_j  \mapsto (1+qX_i) (1+q^3X_i) X_j \\ X'_k \mapsto
(1+qX_i^{-1})^{-1}X_k \quad X'_m \mapsto (1+qX_i^{-1})^{-1}X_m
\end{gather*}

\item[Case 5] (the inverse of Case 4) if $\lambda_k$ is identified
with $\lambda_m$, and $\lambda_j$ is distinct from $\lambda_l$,
then
\begin{gather*}
    X'_j  \mapsto  (1+qX_i)X_j
\quad X'_l
    \mapsto  (1+qX_i)X_l \\
    X'_k \mapsto
(1+qX_i^{-1})^{-1}(1+q^3X_i^{-1})^{-1}X_k
\end{gather*}

\item[Case 6] if $\lambda_j$ is identified with $\lambda_k$, and
$\lambda_l$ is identified with $\lambda_m$ (in which case $S$ is a
$3$--times punctured sphere), then
$$
    X'_j  \mapsto  X_iX_j \quad X'_l  \mapsto
X_iX_l;
$$

\item[Case 7] (the inverse of Case 6) if $\lambda_j$ is identified
with $\lambda_m$, and $\lambda_k$ is identified with $\lambda_l$
(in which case $S$ is a $3$--times punctured sphere), then
$$
    X'_j  \mapsto  X_iX_j
\quad X'_k  \mapsto  X_iX_k;
$$

\item[Case 8] if $\lambda_j$ is identified with $\lambda_l$, and
$\lambda_k$ is identified with $\lambda_m$ (in which case $S$ is a
once punctured torus), then
\begin{gather*}
    X'_j  \mapsto
(1+qX_i)(1+q^3X_i)X_j \\
    X'_k
\mapsto (1+qX_i^{-1})^{-1}(1+q^3X_i^{-1})^{-1}X_k
\end{gather*}

\end{description}

\end{definition}

It turns out that the maps $\widehat{\alpha}$ and $\widehat{\Delta}_i$ can be extended to the whole algebra $\widehat{\mathcal{T}}^q_{\lambda'}$ as algebra homomorphisms from $\widehat{\mathcal{T}}^q_{\lambda'}$ to $\widehat{\mathcal{T}}^q_{\lambda}$.

The motivation of the definition of $\widehat{\alpha}$ and $\widehat{\Delta}_i$ is that they are reduced to the corresponding shear coordinate changes (Proposition \ref{prop:DiagExch}) when $q=1$. 

\begin{prop}[Liu \cite{Liu1}] If an ideal triangulation $\lambda'$ is obtained from another one $\lambda$ by an operation $\pi,$ where $\pi=\alpha$ for some $\alpha\in \mathfrak{S}_{3m},$ or $\pi=\Delta_i$ for some $i$, then $\widehat{\pi}: \widehat{\mathcal{T}}^q_{\lambda'}
\rightarrow \widehat{\mathcal{T}}^q_{\lambda}$ as in Definition \ref{def:exchange} is an isomorphism between the two algebras.
\end{prop}

\begin{prop}[Liu \cite{Liu1}] The map $\widehat{\alpha}$ and $\widehat{\Delta}_i$ satisfy the following relations which correspond to the relations in Lemma \ref{thm:graph-ideal}:

\begin{enumerate}
\item $\widehat{\alpha\beta} = \widehat{\alpha}\circ\widehat{\beta}$ for every $\alpha$, $\beta\in\mathfrak{S}_{3m}$;

\item $\widehat{\Delta}_i\circ\widehat{\Delta}_i=\mathrm{Id}$;

\item $\widehat{\alpha} \circ  \widehat{\Delta}_i = \widehat{\Delta}_{\alpha(i)} \circ  \widehat{\alpha}$ for every $\alpha\in\mathfrak{S}_{3m}$;

\item If $\lambda_i$ and $\lambda_j$ do not belong to the same triangle of $\lambda\in \Lambda(S)$, then $\widehat{\Delta}_i \circ \widehat{\Delta}_j= \widehat{\Delta}_j \circ \widehat{\Delta}_i$;

\item If three triangles of an ideal triangulation $\lambda\in \Lambda(S)$ form a pentagon with
diagonals $\lambda_i$, $\lambda_j$ as in Figure \ref{fig:pentagon-ideal}, then
\begin{equation}
\widehat{\Delta}_i\circ \widehat{\Delta}_j\circ \widehat{\Delta}_i\circ \widehat{\Delta}_j\circ \widehat{\Delta}_i = \widehat{\alpha}_{i\leftrightarrow j}.
\end{equation}

\end{enumerate}

\end{prop}

\subsection{The quantum Teichm\"uller space}

\begin{theorem}[Liu \cite{Liu1}]\label{thm:main-ideal}
There is a family of algebra isomorphisms
$$\Phi_{\lambda\lambda'}^q:
\widehat{\mathcal{T}}^q_{\lambda'} \rightarrow
\widehat{\mathcal{T}}^q_{\lambda}$$ defined as $\lambda$,
$\lambda' \in \Lambda(S)$ ranges over all pairs of ideal
triangulations, such that:
\begin{enumerate}
\item $\Phi_{\lambda\lambda''}^q = \Phi_{\lambda\lambda'}^q \circ
\Phi_{\lambda'\lambda''}^q$ for every $\lambda$, $\lambda'$,
$\lambda''\in \Lambda(S)$;

\item $\Phi_{\lambda\lambda'}^q$ is the isomorphism defined in Definition \ref{def:exchange} when $\lambda'$ is obtained
from $\lambda$ by a reindexing or a diagonal exchange.

\item $\Phi_{\lambda\lambda'}^q$ depends only on $\lambda$ and $\lambda'$.
\end{enumerate}
\end{theorem}

The \emph{quantum (enhanced) Teichm\"uller space}\index{Teichm\"uller space!quantum} of $S$ can now
be defined as the algebra
$$
\widehat {\mathcal{T}}^q_S= \bigg(
\bigsqcup_{\lambda\in\Lambda(S)}
\widehat{\mathcal{T}}^q_{\lambda}\bigg)/\sim
$$
where the relation $\sim$ is defined by the property that, for
$X\in \widehat{\mathcal{T}}^q_{\lambda}$ and $X'\in
\widehat{\mathcal{T}}^q_{\lambda'}$,
$$
X \sim X' \Leftrightarrow X=\Phi^q_{\lambda,\lambda'}(X').
$$

The quantum Teichm\"uller space $\widehat{\mathcal T}_S^q$
is a noncommutative deformation of the algebra of rational functions on the enhanced Teichm\"uller space $\widetilde{\mathcal
T}(S)$.

\section{The quantum Teichm\"uller space: properties }

In this section we survey some interesting properties and applications of the quantum Teichm\"uller space. The uniqueness of the construction of the quantum Teichm\"uller space is established by H. Bai \cite{Bai}. In \cite{BBL,Bon2,BonLiu, Liu2}, it is shown that the quantum Teichm\"uller space $\widehat{\mathcal T}_S^q$ has a rich representation theory which also produces an invariant of hyperbolic 3-manifolds. 

We would like to mention the following related important works without providing more details. 

H. Bai \cite{Bai2} shows that Kashaev's $6j$-symbols \cite{K-knot1, K-knot2} are intertwining operators of local representations of quantum Teichm\"ller spaces introduced in \cite{BBL}. Note that appearance of Kashaev's $6j$-symbols in quantum Teichm\"ller theory at roots of unity is already explicit in \cite{Kas1}(see the operator $T_{h,x,y}$ in Proposition 10).

C. Hiatt \cite{Hia} proves that for the torus with one hole and $p \geq 1$ punctures and the sphere with four holes
there is a family of quantum trace functions in the quantum Teichm\"uller space, analog to the non-quantum
trace functions in Teichm\"uller space, satisfying the properties proposed by Chekhov and Fock in \cite{CF}.

For a punctured surface $S$, a point of its Teichm\"uller space
${\mathcal T}(S)$ determines an irreducible representation of its quantization ${\mathcal T}_S^q$ .
J. Roger \cite{R} analyzes the behavior of these representations as one goes to infinity
in ${\mathcal T}_S$. He shows that an irreducible representation of ${\mathcal T}_S^q$ limits to
a direct sum of representations of ${\mathcal T}_{S_\gamma}^q$, where $S_\gamma$
is obtained from $S$ by pinching a multicurve 
$\gamma$ to a set of nodes. The result is analogous to
the factorization rule found in conformal field theory.

The skein algebra and the quantum Teichm\"uller space are considered as two different quantizations of the character variety
consisting of all representations of surface groups in $\mbox{PSL}_2(\mathbb C)$.
F. Bonahon and H. Wong \cite{BW1,BW2} construct a homomorphism from the skein
algebra to the quantum Teichm\"uller space which, when restricted the classical
case, corresponds to the equivalence between these two algebras through trace
functions.

\subsection{Uniqueness}

The original definition of the quantum
Teichm\"uller space was motivated by geometry. However, H. Bai \cite{Bai}
shows that it is intrinsically tied to the combinatorics of the
set $\Lambda(S)$. Indeed, H. Bai proves that the coordinate change isomorphisms
considered in Definition \ref{def:exchange} are the only ones which satisfy
a certain number of natural conditions. 

The \textit{discrepancy span} $D(\lambda,
\lambda')$ of two ideal triangulations $\lambda$, $\lambda'$ is
the closure of the union of those connected components of
$S-\lambda$ which are not isotopic to a component of $S-\lambda'$.

The coordinate change
isomorphisms $\Phi_{\lambda \lambda'}^q$ are said to satisfy the
\textit{Locality Condition} if the following holds. 
Let $\lambda$ and
$\lambda'$ be two ideal triangulations indexed in such a way that
$\lambda_i \subset D(\lambda,
\lambda')$ when  $i \leq k$, and $\lambda_i' = \lambda_i$ when
$i>k$. Then
\begin{enumerate}
   \item $\Phi_{\lambda
\lambda'}^q(X_i')=X_i$ for every  $i>k$;
   \item $\Phi_{\lambda
\lambda'}^q(X_i')=f_i(X_1,X_2,\cdots,X_k)$ for every $i\leq k$,
where $f_i$ is a multi-variable rational function depending only on the
combinatorics of $\lambda$ and $\lambda'$ in $D(\lambda, \lambda')$ in
the following sense: For any two pairs of ideal triangulation
   $(\lambda, \lambda')$, $(\lambda'', \lambda''')$ for which
there exists a diffeomorphism $\psi: D(\lambda,
\lambda')\rightarrow D(\lambda'', \lambda''')$ sending $\lambda_i$
to $\lambda_j''$ and $\lambda_j'$ to $\lambda_j'''$ for every
$1\leq j\leq k$, then
\begin{align*}
\Phi_{\lambda
\lambda'}^q(X_i') &=f_i(X_1,X_2,\cdots,X_k) \textrm{ and}\\
\Phi_{\lambda''
\lambda'''}^q(X_i''') &=f_i(X_1'',X_2'',\cdots,X_k'')
\end{align*}
for the same
rational function $f_i$.
\end{enumerate}

\begin{prop}[Bai \cite{Bai}] The algebra isomorphisms 
$\Phi_{\lambda\lambda'}^q: \widehat{\mathcal{T}}^q_{\lambda'} \rightarrow \widehat{\mathcal{T}}^q_{\lambda}$ in Theorem \ref{thm:main-ideal} satisfies the Locality Condition.
\end{prop}

\begin{theorem}[Bai \cite{Bai}]
\label{thm:CoordChangesUnique}
Assume that the surface $S$ satisfies $\chi(S)<-2$. Then the family of coordinate change isomorphisms
$\Phi_{\lambda\lambda'}^q$ in Theorem \ref{thm:main-ideal} is unique up to a uniform rescaling
and/or inversion of the $X_i$.

Namely, if
$$\Psi_{\lambda
\lambda'}^q:\mathbb{C}(X_1',X_2',\dots,X_n')_{\lambda'}^q
\rightarrow \mathbb{C}(X_1,X_2,\dots,X_n)_{\lambda}^q$$ 
is another family of
isomorphisms satisfying the conditions of
Theorem~\ref{thm:main-ideal} and the Locality Condition, then there exists a non-zero
constant $\xi \in \mathbb{C}(q)$ and a sign $\varepsilon =\pm1$
such that $\Psi_{\lambda \lambda'}^q=\Theta_{\lambda} \circ
\Phi_{\lambda \lambda'}^q \circ \Theta_{\lambda'}^{-1}$ for any
pair of ideal triangulations $\lambda$, $\lambda'$, where
$\Theta_{\lambda}: \mathbb{C}(X_1,X_2,\dots,X_n)_{\lambda}^q
\rightarrow \mathbb{C}(X_1,X_2,\dots,X_n)_{\lambda}^q$ is the
isomorphism defined by $\Theta_{\lambda}(X_i)=\xi X_i^\varepsilon$
for every $i$.
\end{theorem}

Theorem~\ref{thm:CoordChangesUnique} is false when $S$ is the
once-punctured torus or the 3--times punctured sphere. The
uniqueness property for the twice-punctured torus and
the 4--times punctured sphere has not been established.

\subsection{Representations}

In this subsecton, our exposition follows \cite{BonLiu} closely.

A standard method to move from abstract algebraic constructions to more concrete
applications is to consider finite-dimensional representations. In the case
of algebras, this means algebra homomorphisms valued in the algebra $\mbox{End}(V)$ of
endomorphisms of a finite-dimensional vector space $V$ over $\mathbb{C}$. Elementary considerations
show that these can exist only when $q$ is a root of unity.

\begin{theorem}[Bonahon-Liu \cite{BonLiu}]
\label{thm:RepCheFockIntro}
Suppose that $q^2$ is a primitive $N$--th root of unity. For any ideal triangulation $\lambda$ of a surface $S$, every irreducible finite-dimensional representation of the Chekhov-Fock algebra $ \mathcal T^q_\lambda$ has dimension $N^{3g+p-3}$ if $N$ is odd, and
$N^{3g+p-3}/2^g$ if $N$ is even. Up to isomorphism, such a representation is classified by:
\begin{enumerate}
\item a non-zero complex number $x_i \in \mathbb C^*$ associated to each edge of $\lambda$;

\item a choice of an $N$--th root for each of $p$ explicit monomials in the numbers $x_i$;

\item when $N$ is even, a choice of square root for each of $2g$ explicit monomials in the numbers $x_i$. 
\end{enumerate} 
Conversely, any such data can be realized by an irreducible finite-dimensional representation of $ \mathcal T^q_\lambda$.
\end{theorem}

Theorem~\ref{thm:RepCheFockIntro} shows that the Chekhov-Fock algebra has a rich representation theory. Unfortunately, for dimension reasons, its fraction algebra $\widehat{\mathcal T}^q_\lambda$ and, consequently, the quantum Teichm\"uller space $\widehat{\mathcal T}^q_S$ cannot have any finite-dimensional representation. This leads us to introduce the \emph{polynomial core}\index{polynomial core} $\mathcal T_S^q$ of the quantum Teichm\"uller space
$\widehat{\mathcal T}^q_S$, defined as the family $\{ \mathcal T_\lambda^q\} _{\lambda \in \Lambda(S)}$ of all Chekhov-Fock algebras $\mathcal T_\lambda^q$, considered as subalgebras of $\widehat{\mathcal T}^q_S$, as $\lambda$ ranges over the set $\Lambda(S)$ of all isotopy classes of ideal triangulations of the surface $S$.  

Theorem~\ref{thm:RepCheFockIntro} says that, up to a finite number of choices,
an irreducible representation of 
$\mathcal T^q_\lambda$ is
classified by certain numbers
$x_i\in \mathbb C^*$ associated
to the edges of the ideal
triangulation 
$\lambda$ of $S$. There is a
classical geometric object which is
also associated to $\lambda$ with
the same edge weights $x_i$.
Namely, we can consider in the hyperbolic
3--space
$\mathbb H^3$ the
pleated surface that has pleating
locus
$\lambda$, that has shear
parameter along the
$i$--th edge of $\lambda$ equal to
the real part of $\ln x_i$, and
that has bending angle along this
edge equal to the imaginary part
of $\ln x_i$. In turn, this
pleated surface has a
\emph{monodromy representation}\index{monodromy representation},
namely a group homomorphism from
the fundamental group $\pi_1(S)$
to the  group
$\mathrm{Isom}^+(\mathbb H^3)
\cong \mathrm{PSL}_2(\mathbb C)$
of orientation-preserving
isometries of $\mathbb H^3$. This
construction associates to a
representation of the Chekhov-Fock
algebra $ \mathcal T^q_\lambda$ a
group homomorphism
$r: \pi_1(S) \rightarrow
\mathrm{PSL}_2(\mathbb C)$,
well-defined up to conjugation by
an element of
$\mathrm{PSL}_2(\mathbb C)$.

\begin{theorem} [Bonahon-Liu \cite{BonLiu}]
\label{thm:HyperShadowIntro} Let
$q$ be a primitive $N$--th root
of $(-1)^{N+1}$, for instance
$q=-\mathrm e^{2 \pi
\mathrm i/N}$. If 
$\rho =
\{\rho_\lambda:
\mathcal T^q_\lambda
\rightarrow
\mathrm{End} (V)\}_{\lambda \in
\Lambda(S)}$ is a
finite-dimensional irreducible
representation of the polynomial
core $\mathcal T^q_S$ of the
quantum Teichm\"uller space
$\widehat{\mathcal T}^q_S$, the
representations $\rho_\lambda$
induce the same monodromy
homomorphism $r_\rho: \pi_1(S)
\rightarrow
\mathrm{PSL}_2(\mathbb C)$.
\end{theorem}

Theorem~\ref{thm:HyperShadowIntro}
is essentially equivalent to the
property that, for the choice of
$q$ indicated, the pleated
surfaces respectively associated
to the representations 
$\rho_\lambda:
\mathcal T^q_\lambda
\rightarrow
\mathrm{End} (V)$ and 
$\rho_\lambda \circ
\Phi_{\lambda\lambda'}^q:
\mathcal T^q_{\lambda'}
\rightarrow
\mathrm{End} (V)$ have (different
pleating loci but) the same
monodromy representation 
$r_\rho: \pi_1(S)
\rightarrow
\mathrm{PSL}_2(\mathbb C)$.

The homomorphism $r_\rho$ is the
\emph{hyperbolic shadow}\index{hyperbolic shadow} of the
representation $\rho$. Not every
homomorphism $r: \pi_1(S)
\rightarrow
\mathrm{PSL}_2(\mathbb C)$ is the
hyperbolic shadow of a
representation of the polynomial
core, but many of them are:

\begin{theorem}[Bonahon-Liu \cite{BonLiu}]
\label{thm:HyperGivesRepIntro} An
injective homomorphism  
$r: \pi_1(S)
\rightarrow
\mathrm{PSL}_2(\mathbb C)$ is the
hyperbolic shadow of a finite
number of irreducible
finite-dimensional
representations of the polynomial
core $\mathcal T^q_S$, up to
isomorphism. More precisely, this
number of representations is
equal to $2^lN^p$ if $N$ is odd,
and $ 2^{2g+l}N^p$ if $N$ is
even, where $l$ is the number
of ends of $S$ whose image under
$r$ is loxodromic.
\end{theorem}

Let $\varphi$ be a diffeomorphism
of the surface $S$. Suppose in
addition that $\varphi$ is
homotopically aperiodic
(also called homotopically
pseudo-Anosov), so that its
(3--dimensional) mapping torus
$M_\varphi$ admits a complete
hyperbolic metric. The hyperbolic
metric of $M_\varphi$ gives an
injective homomorphism 
$r_\varphi: \pi_1(S)
\rightarrow
\mathrm{PSL}_2(\mathbb C)$ such
that $r_\varphi \circ \varphi^*$
is conjugate to $r_\varphi$, where
$\varphi^*$ is the isomorphism of
$\pi_1(S)$ induced by $\varphi$.

The diffeomorphism $\varphi$ also
acts on the quantum Teichm\"uller
space and on its polynomial core $
\mathcal T_S^q$. In particular,
it acts on the set of
representations of
$\mathcal T_S^q$ and, because
$r_\varphi
\circ
\varphi^*$ is conjugate to
$r_\varphi$, it sends a
representation with hyperbolic
shadow $r_\varphi$ to another
representation with shadow
$r_\varphi$. 
Actually, when $N$ is odd,
there is a preferred
representation 
$\rho_\varphi$ of $\mathcal
T_S^q$ which
is fixed by the action of
$\varphi$, up to isomorphism. 
This statement means that, for
every ideal triangulation
$\lambda$, we have a
representation
$\rho_\lambda:
\mathcal T^q_\lambda
\rightarrow
\mathrm{End} (V)$ of dimension 
$N^{3g+p-3}$ and an isomorphism
$L_\varphi^q$ of $V$ such that
\begin{equation*}
\rho_{\varphi(\lambda)} \circ
\Phi_{\varphi(\lambda) \lambda}
(X) = L_\varphi^q \cdot
\rho_\lambda (X)
\cdot (L_\varphi^q)^{-1}
\end{equation*} in
$\mathrm{End}(V)$ for every
$X\in \mathcal T_\lambda^q$, for
a suitable interpretation of the
left hand side of the equation.

\begin{theorem} [Bonahon-Liu \cite{BonLiu}]
\label{thm:intro:Invariant}
Let $N$ be odd. Up to
conjugation and up to
multiplication by a constant, the
isomorphism
$L_\varphi^q$ depends only on the
homotopically aperiodic
diffeomorphism 
$\varphi: S \rightarrow S$ and
on the primitive
$N$--th root $q$ of
$1$. 
\end{theorem}

Explicit computations of these invariants for the once-punctured torus or the 4--times punctured sphere
are provided in X. Liu \cite{Liu2}. 

H. Bai, F. Bonahon and X. Liu \cite{BBL} investigate another type of representations of the quantum Teichm\"uller space, called local representations, which are somewhat simpler to analyze and more closely connected to the combinatorics of ideal triangulations.

\section{Kashaev algebra}

In this section, we establish the theory of Kashaev algebra\index{Kashaev algebra} which is parallel to the theory of the quantum Teichm\"uller space. The exposition follows \cite{GL} closely.

\subsection{Decorated ideal triangulations}

Let $S$ be an oriented surface of genus $g$ with $p\geq 1$ punctures and negative Euler characteristic, i.e., $m=2g-2+p>0.$  

A decorated ideal triangulation\index{ideal triangulation!decorated} $\tau$ of $S$ introduced by Kashaev \cite{Kas1} is an ideal triangulation such that the ideal triangles are numerated as $\{\tau_1,\tau_2,...,\tau_{2m}\}$ and there is a mark (a star symbol) at a corner of each ideal triangle. Denote by $\triangle(S)$ the set of isotopy classes of decorated ideal triangulations of the surface $S$.

The set $\triangle(S)$ admits a natural action of the symmetric group on the set $\{1,2,...,2m\}$,
$\mathfrak{S}_{2m}$, acting by
permuting the indexes of the ideal triangles of $\tau$. Namely $\tau'
= \alpha(\tau)$ for $\alpha \in\mathfrak{S}_{2m}$ if  $\tau_i=\tau'_{\alpha(i)}$.

\begin{figure}[ht!]
\labellist\small\hair 2pt
\pinlabel $(1)$ at 205 200
\pinlabel $\longrightarrow$ at 205 81
\pinlabel $\varphi_{ij}$ at 205 99
\pinlabel $\tau_i$ at 50 85
\pinlabel $\tau_j$ at 107 85
\pinlabel $*$ at 13 85
\pinlabel $*$ at 144 85
\pinlabel $\tau'_i$ at 329 107
\pinlabel $\tau'_j$ at 329 48
\pinlabel $*$ at 329 144
\pinlabel $*$ at 329 13
\pinlabel $(2)$ at 205 0
\pinlabel $\longrightarrow$ at 205 273
\pinlabel $\rho_i$ at 205 291
\pinlabel $\tau_i$ at 80 264
\pinlabel $\tau'_i$ at 329 264
\pinlabel $*$ at 79 301
\pinlabel $*$ at 260 240

\endlabellist
\centering
\includegraphics[scale=0.5]{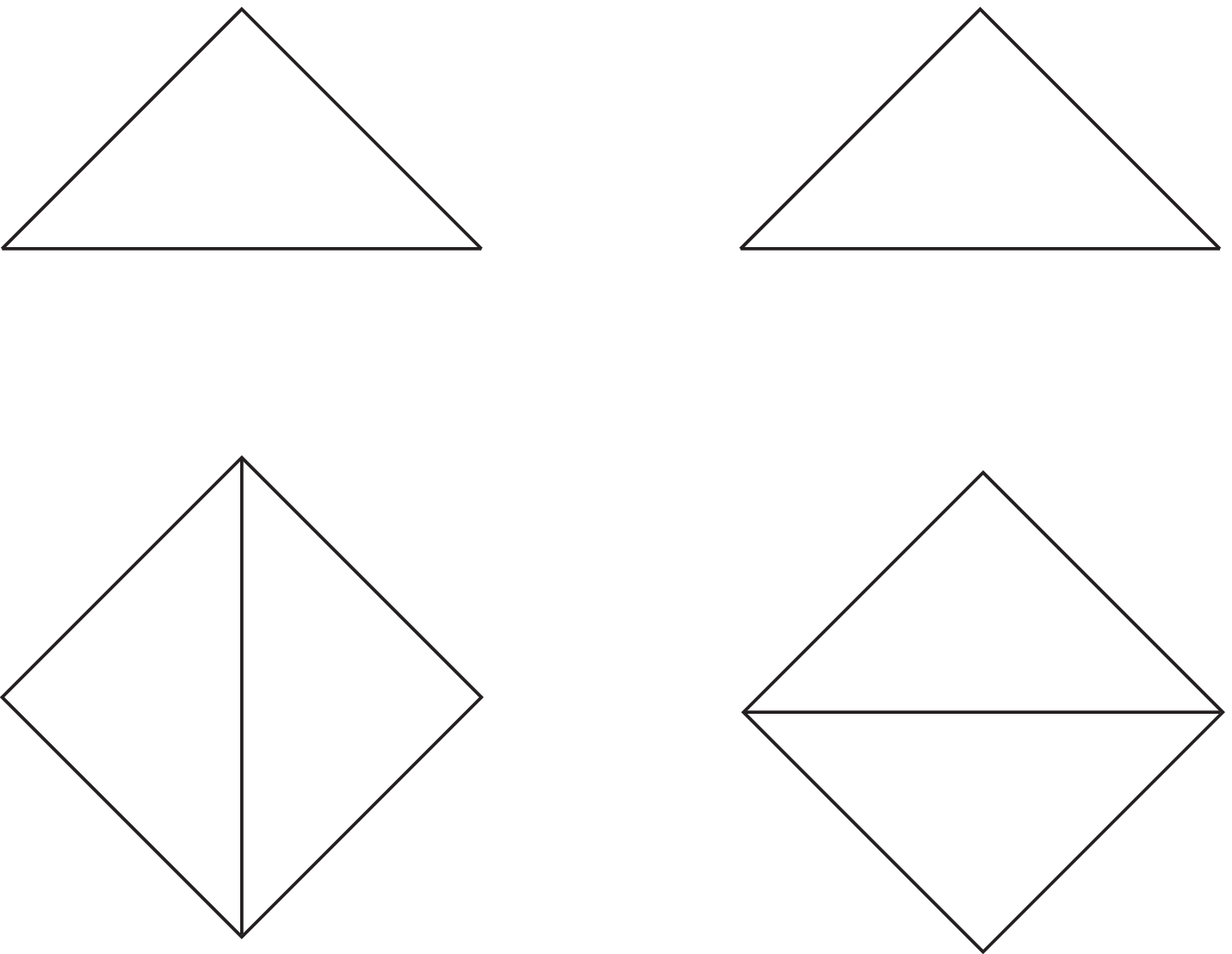}
\caption{}
\label{fig:maps}
\end{figure}

Another important transformation of $\triangle(S)$ is provided by
the \emph {diagonal exchange}\index{diagonal exchange} $\varphi_{ij}:
\triangle(S)\rightarrow \triangle(S)$ defined as follows. Suppose that two ideal triangles $\tau_i, \tau_j$ share an edge $e$ such that the marked corners are opposite to the edge $e.$ Then $\varphi_{ij}(\tau)$ is obtained by rotating the interior of the union $\tau_i\cup \tau_j$ $90^\circ$ in the clockwise order, as illustrated in Figure \ref{fig:maps}(2).

The last transformation of $\triangle(S)$ is the \emph {mark rotation}\index{mark rotation} $\rho_i: \triangle(S)\rightarrow \triangle(S)$. $\rho_i(\tau)$ is obtained by relocating the mark of the ideal triangle $\tau_i$ from one corner to the next corner in the counterclockwise order, as illustrated in Figure \ref{fig:maps}(1).

\begin{figure}[ht!]
\labellist\small\hair 2pt
\pinlabel $\tau_i$ at 33 81
\pinlabel $\tau_j$ at 81 42
\pinlabel $\tau_k$ at 127 81
\pinlabel $*$ at 7 92
\pinlabel $*$ at 42 7
\pinlabel $*$ at 133 17

\endlabellist
\centering
\includegraphics[scale=0.5]{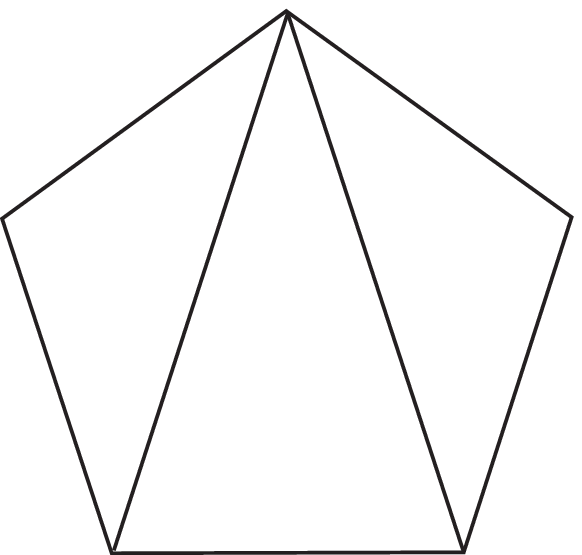}
\caption{}
\label{fig:pentagon}
\end{figure}

\begin{lemma}\label{thm:graph}
The reindexings, diagonal exchanges and mark rotations satisfy the following
relations:
\begin{enumerate}

\item $(\alpha\beta)(\tau) = \alpha(\beta(\tau))$ for every $\alpha$, $\beta\in\mathfrak{S}_{2m}$;

\item $\varphi_{ij}\circ \varphi_{ij}=\alpha_{i\leftrightarrow j},$ where $\alpha_{i\leftrightarrow j}$ denotes the transposition exchanging $i$ and $j$;

\item $\alpha \circ \varphi_{ij} = \varphi_{\alpha(i)\alpha(j)}\circ \alpha$ for every $\alpha\in\mathfrak{S}_{2m}$;

\item $\varphi_{ij}\circ\varphi_{kl}(\tau)=\varphi_{kl}\circ \varphi_{ij}(\tau)$, for $\{i,j\}\neq\{k,l\}$;

\item If three triangles $\tau_i,\tau_j,\tau_k$ of an ideal
triangulation $\tau \in \triangle(S)$ form a pentagon and their marked corners are located as in Figure \ref{fig:pentagon}, then the Pentagon Relation holds:
$$\omega_{jk}\circ\omega_{ik}\circ\omega_{ij}(\tau)=\omega_{ij}\circ\omega_{jk}(\tau),$$
where $\omega_{\mu\nu}=\rho_\mu\circ\varphi_{\mu\nu}\circ\rho_\nu$;

\item $\rho_i \circ \rho_i \circ \rho_i=\mathrm{Id}$;

\item $\rho_i \circ \rho_j= \rho_j \circ \rho_i$;

\item $\alpha \circ \rho_i= \rho_{\alpha(i)} \circ \alpha$ for for every $\alpha\in\mathfrak{S}_{2m}$.

\end{enumerate}

\end{lemma}

Lemma \ref{thm:graph} is essentially contained in Kashaev \cite{Kas2} where $\omega_{ij}$ is used as the diagonal exchange.

The following two results about decorated ideal triangulations can be easily proved using Penner's result about ideal triangulations \cite{Pen}.

\begin{theorem}\label{thm:Penner1} Given two decorated ideal triangulations
$\tau,\tau'\in \triangle(S)$, there exists a finite sequence
of decorated ideal triangulations $\tau=\tau_{(0)}$, $\tau_{(1)}$,
\dots, $\tau_{(n)}=\tau'$ such that each $\tau_{(k+1)}$
is obtained from $\tau_{(k)}$ by a diagonal exchange or by a mark rotation or by a
reindexing of its ideal triangles.
\end{theorem}

\begin{theorem}\label{thm:Penner2} Given two decorated ideal triangulations
$\tau,\tau'\in \triangle(S)$ and given two sequences
$\tau=\tau_{(0)}$, $\tau_{(1)}$, \ldots, $\tau_{(n)}
=\tau'$ and $\tau=\overline{\tau}_{(0)}$, $\overline{\tau}_{(1)}$, \ldots,
$\overline{\tau}_{(\overline{n})} =\tau'$ of diagonal exchanges, mark rotations and reindexings
connecting them as in Theorem \ref{thm:Penner1}, these two
sequences can be related to each other by successive applications
of the following moves and of their inverses. These moves correspond to the relations in Lemma \ref{thm:graph}.
\begin{enumerate}

\item Replace \dots,
$\tau_{(k)}$, $\beta(\tau_{(k)})$, $\alpha\circ \beta(\tau_{(k)})$, \dots

\ \ \ by \dots, $\tau_{(k)}$, $(\alpha \beta)(\tau_{(k)})$, \dots
where $\alpha$, $\beta\in \mathfrak{S}_n$.

\item Replace \dots,
$\tau_{(k)}$, $\varphi_{ij}(\tau_{(k)}) $, $\varphi_{ij}\circ \varphi_{ij}(\tau_{(k)}) $ \dots

\ \ \ by \dots, $ \tau_{(k)}$, $\alpha_{i\leftrightarrow j}(\tau_{(k)})$, \dots .

\item Replace \dots,
$\tau_{(k)}$, $\varphi_{ij}(\tau_{(k)}) $, $\alpha \circ \varphi_{ij}(\tau_{(k)})$, \dots

\ \ \ by \dots, $\tau_{(k)}$, $\alpha(\tau_{(k)})$, $\varphi_{\alpha(i)\alpha(j)}\circ \alpha(\tau_{(k)})$, \dots
where $\alpha\in \mathfrak{S}_n$.

\item Replace \dots, $\tau_{(k)}$, $\varphi_{kl}(\tau_{(k)})$, $\varphi_{ij}\circ \varphi_{kl}(\tau_{(k)})$, \dots

\ \ \ by \dots, $ \tau_{(k)}$, $\varphi_{ij}(\tau_{(k)})$, $\varphi_{kl}\circ \varphi_{ij}(\tau_{(k)})$, \dots  where $\{i,j\}\neq\{k,l\}$.

\item Replace \dots, $\tau_{(k)}$, $\omega_{ij}(\tau_{(k)})$, $\omega_{ik} \circ\omega_{ij}(\tau_{(k)})$,
$\omega_{jk} \circ\omega_{ik} \circ\omega_{ij}(\tau_{(k)})$, \dots,

\ \ \ by \dots, $\tau_{(k)}$, $\omega_{jk}(\tau_{(k)})$, $\omega_{ij} \circ\omega_{jk}(\tau_{(k)})$, \dots
where $\omega_{\mu\nu}=\rho_\mu\circ\varphi_{\mu\nu}\circ\rho_\nu$.

\item Replace \dots, $\tau_{(k)}$,  $\rho_i(\tau_{(k)})$,  $\rho_i\circ \rho_i(\tau_{(k)})$, $\tau_{(k)}$ \dots

\ \ \ by  \dots, $\tau_{(k)}$, \dots.

\item Replace \dots, $\tau_{(k)}$, $\rho_i(\tau_{(k)})$, $\rho_j\circ \rho_i(\tau_{(k)})$, \dots

\ \ \ by \dots, $\tau_{(k)}$, $\rho_j(\tau_{(k)})$, $\rho_i\circ \rho_j(\tau_{(k)})$, \dots.

\item Replace \dots, $\tau_{(k)}$, $\rho_i(\tau_{(k)})$, $\alpha\circ \rho_i(\tau_{(k)})$, \dots

\ \ \ by \dots, $\tau_{(k)}$, $\alpha(\tau_{(k)})$, $\rho_{\alpha(i)}\circ \alpha(\tau_{(k)})$, \dots.

\end{enumerate}
\end{theorem}

\subsection{Kashaev coordinates}

For a decorated ideal triangulation $\tau$ of a punctured surface $S$, Kashaev \cite{Kas1} associated to each ideal triangle $\tau_i$ two numbers $\ln y_i, \ln z_i$. A Kashaev coordinate\index{Kashaev coordinate} is a vector $(\ln y_1, \ln z_1,...,\ln y_{2m}, \ln z_{2m})\in \mathbb{R}^{4m}.$

Denote by $(y_1,z_1,...,y_{2m},z_{2m})$ the exponential Kashaev coordinate for the decorated ideal triangulation $\tau$. Denote by $(y'_1,z'_1,...,y'_{2m},z'_{2m})$ the exponential Kashaev coordinate for the decorated ideal triangulation $\tau'$. Kashaev \cite{Kas1} introduces the change of coordinates as follows.

\begin{definition}[Kashaev \cite{Kas1}]\label{def:coor-change}
Suppose that a decorated ideal triangulation $\tau'$ is obtained from another one $\tau$ by reindexing the ideal triangles, i.e., $\tau'=\alpha(\tau)$ for some $\alpha\in \mathfrak{S}_{2m},$ then we define $(y'_i, z'_i)=(y_{\alpha(i)}, z_{\alpha(i)})$ for any $i=1,...,2m.$

Suppose that a decorated ideal triangulation $\tau'$ is obtained from another one $\tau$ by a mark rotation, i.e., $\tau'=\rho_i(\tau)$ for some $i,$ then we define $(y'_j, z'_j)=(y_j,z_j)$ for any $j\neq i$ while
$$(y'_i, z'_i)=(\frac{z_i}{y_i}, \frac1{y_i}).$$

Suppose a decorated ideal triangulation $\tau'$ is obtained from another one $\tau$ by a diagonal exchange, i.e., $\tau'=\varphi_{ij}(\tau)$ for some $i,j,$ then we define $(y'_k, z'_k)=(y_k,z_k)$ for any $k\notin \{i,j\}$ while
$$(y'_i,z'_i,y'_j,z'_j)=
(\frac{z_j}{y_iy_j+z_iz_j},\frac{y_i}{y_iy_j+z_iz_j},\frac{z_i}{y_iy_j+z_iz_j},\frac{y_j}{y_iy_j+z_iz_j}).$$
\end{definition}

Kashaev \cite{Kas1} considered $\omega_{ij}$ instead of $\varphi_{ij}$.

There is a natural relationship between Kashaev coordinates and Penner coordinates\index{Penner coordinate} for the decorated Teichm\"uller space\index{Teichm\"uller space!decorated} which is established in \cite{Kas1}. For an exposition, see also Teschner \cite{Te}. 

\subsection{Generalized Kashaev algebra: triangulation-dependent}

For a decorated ideal triangulation $\tau$ of a punctured surface $S$, Kashaev \cite{Kas1} introduced an algebra $\mathcal K^q_{\tau}$ on $\mathbb C$ generated by $Y_1^{\pm},Z_1^{\pm},Y_2^{\pm},Z_2^{\pm},...,Y_{2m}^{\pm},Z_{2m}^{\pm},$ with $Y_i^{\pm},Z_i^{\pm}$ associated to an ideal triangle $\tau_i,$ and by the relations:
\begin{equation*}
\begin{split}
Y_iY_j&=Y_jY_i, \\
Z_iZ_j&=Z_jZ_i, \\
Y_iZ_j&=Z_jY_i\ \ \mbox{if}\ \ i\neq j, \\
Z_iY_i&=q^2Y_iZ_i.
\end{split}
\end{equation*}

Kashaev's original definition is $Y_iZ_i=q^2Z_iY_i.$ We adopt our convention to make it compatible with the quantum Teichm\"uller space \cite{Liu1}. Kashaev's parameter $q$ corresponds to our $q^{-1}$.

The algebra $\widehat{\mathcal K}^q_{\tau}$ is the fraction division algebra of $\mathcal K^q_{\tau}$.

In particular, when $q=1,$ $\mathcal K^q_{\tau}$ and $\widehat{\mathcal K}^q_{\tau}$ respectively coincide with the Laurent polynomial algebra $\mathbb C [Y_1^{\pm},Z_1^{\pm},...,Y_{2m}^{\pm},Z_{2m}^{\pm}]$ and the rational
fraction algebra $\mathbb C (Y_1,Z_1,...,Y_{2m},Z_{2m})$. The general
$\mathcal K^q_{\tau}$ and $\widehat{\mathcal K}^q_{\tau}$ can be considered as deformations of $\mathcal K^1_{\tau}$ and $\widehat{\mathcal K}^1_{\tau}$.

The algebra $\widehat{\mathcal K}^q_{\tau}$ depends on the decorated ideal triangulation $\tau$. We introduce algebra isomorphisms in the following.

\begin{definition}\label{def:iso} Let $a,b$ be two arbitrary nonzero complex numbers.

Suppose that a decorated ideal triangulation $\tau'$ is obtained from another one $\tau$ by reindexing the ideal triangles, i.e., $\tau'=\alpha(\tau)$ for some $\alpha\in \mathfrak{S}_{2m},$ then we define a map $\widehat{\alpha}$ on the set of generators of $\widehat{\mathcal K}^q_{\tau'}$ to $\widehat{\mathcal K}^q_{\tau}$ by
\begin{align*}
\widehat{\alpha}(Y'_i) &= Y_{\alpha(i)}, \ \ \mbox{for any}\ \ i=1,...,2m,\\
\widehat{\alpha}(Z'_i) &= Z_{\alpha(i)}, \ \ \mbox{for any}\ \ i=1,...,2m.
\end{align*}

Suppose that a decorated ideal triangulation $\tau'$ is obtained from another one $\tau$ by a mark rotation, i.e., $\tau'=\rho_i(\tau)$ for some $i,$ then we define a map $\widehat{\rho}_i$ on the set of generators of $\widehat{\mathcal K}^q_{\tau'}$ to $\widehat{\mathcal K}^q_{\tau}$ by
\begin{align*}
\widehat{\rho}_i(Y'_j) &= Y_j, \ \ \mbox{if}\ \ j\neq i, \\
\widehat{\rho}_i(Z'_j) &= Z_j, \ \ \mbox{if}\ \ j\neq i, \\
\widehat{\rho}_i(Y'_i) &= aY_i^{-1}Z_i, \\
\widehat{\rho}_i(Z'_i) &= Y_i^{-1}.
\end{align*}

Suppose a decorated ideal triangulation $\tau'$ is obtained from another one $\tau$ by a diagonal exchange, i.e., $\tau'=\varphi_{ij}(\tau)$ for some $i,j,$ then we define a map $\widehat{\varphi}_{ij}$ on the set of generators of $\widehat{\mathcal K}^q_{\tau'}$ to $\widehat{\mathcal K}^q_{\tau}$ by
\begin{align*}
\widehat{\varphi}_{ij}(Y_i')&= \ (b Y_iY_j+Z_iZ_j)^{-1}Z_j, \\
\widehat{\varphi}_{ij}(Z_i')&=b(b Y_iY_j+Z_iZ_j)^{-1}Y_i, \\
\widehat{\varphi}_{ij}(Y_j')&= \ (b Y_iY_j+Z_iZ_j)^{-1}Z_i, \\
\widehat{\varphi}_{ij}(Z_j')&=b(b Y_iY_j+Z_iZ_j)^{-1}Y_j.
\end{align*}

\end{definition}

It turns out that the maps $\widehat{\alpha},$ $\widehat{\rho}_i$ and $\widehat{\varphi}_{ij}$ can be extended to the whole algebra $\widehat{\mathcal K}^q_{\tau'}$ as algebra homomorphisms between from $\widehat{\mathcal K}^q_{\tau'}$ to $\widehat{\mathcal K}^q_{\tau}.$

Kashaev \cite{Kas1} considered a special case of these maps when $a=q^{-1}, b=q.$

From the definition, when $q=1,$ we get the coordinate change formula in Definition \ref{def:coor-change}.

\begin{prop}[Guo-Liu \cite{GL}]\label{thm:iso} If a decorated ideal triangulation $\tau'$ is obtained from another one $\tau$ by an operation $\pi,$ where $\pi=\alpha$ for some $\alpha\in \mathfrak{S}_{2m},$ or $\pi=\rho_i$ for some $i$, or $\pi=\varphi_{ij}$ for some $i,j$, then $\widehat{\pi}: \widehat{\mathcal K}^q_{\tau'}\to \widehat{\mathcal K}^q_{\tau}$ as in Definition \ref{def:iso} is an isomorphism between the two algebras.
\end{prop}

\begin{figure}[ht!]
\labellist\small\hair 2pt
\pinlabel $\tau_i$ at 33 81
\pinlabel $\tau_j$ at 81 42
\pinlabel $\tau_k$ at 127 81
\pinlabel $*$ at 7 92
\pinlabel $*$ at 42 7
\pinlabel $*$ at 133 17

\endlabellist
\centering
\includegraphics[scale=0.5]{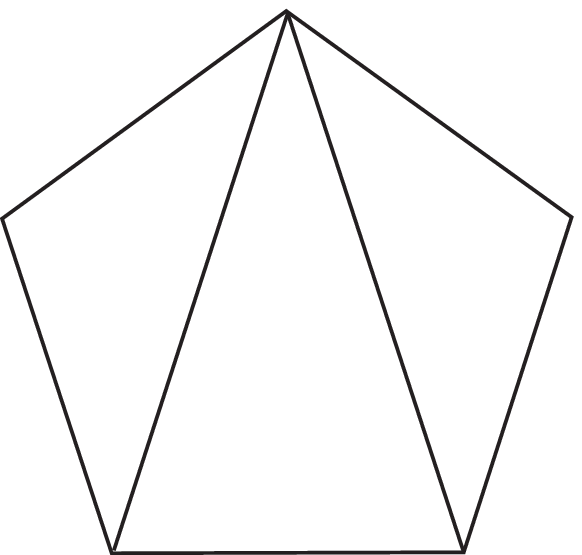}
\caption{(Same as Figure \ref{fig:pentagon})}
\label{fig:pentagon2}
\end{figure}

\begin{prop}[Guo-Liu \cite{GL}]\label{thm:relation} The maps $\widehat{\alpha}, \widehat{\rho}_i$ and $\widehat{\varphi}_{ij}$ satisfy the following relations which correspond to the relations in Lemma \ref{thm:graph}:

\begin{enumerate}

\item $\widehat{\alpha\beta}= \widehat{\alpha}\circ\widehat{\beta}$ for every $\alpha$, $\beta\in\mathfrak{S}_{2m}$;

\item $\widehat{\varphi}_{ij}\circ\widehat{\varphi}_{ij}=\widehat{\alpha}_{i\leftrightarrow j}$;

\item $ \widehat{\alpha} \circ \widehat{\varphi}_{ij}  =
 \widehat{\varphi}_{\alpha(i)\alpha(j)} \circ \widehat{\alpha} $ for every
$\alpha\in\mathfrak{S}_{2m}$;

\item $\widehat{\varphi}_{ij}\circ\widehat{\varphi}_{kl}=\widehat{\varphi}_{kl}\circ \widehat{\varphi}_{ij}$ for $\{i,j\}\neq\{k,l\}$;

\item If three triangles $\tau_i,\tau_j,\tau_k$ of an ideal
triangulation $\tau \in \triangle(S)$ form a pentagon and their marked corners are located as in Figure \ref{fig:pentagon2}, then the Pentagon Relation holds:
$$\widehat{\omega}_{jk}\circ\widehat{\omega}_{ik}\circ\widehat{\omega}_{ij}=
\widehat{\omega}_{ij}\circ\widehat{\omega}_{jk},$$
where $\widehat{\omega}_{\mu\nu}=\widehat{\rho}_\mu\circ\widehat{\varphi}_{\mu\nu}\circ\widehat{\rho}_\nu$;

\item $\widehat{\rho}_i\circ\widehat{\rho}_i\circ\widehat{\rho}_i=\mathrm{Id}$;

\item $\widehat{\rho}_i\circ\widehat{\rho}_j=\widehat{\rho}_j\circ\widehat{\rho}_i$;

\item $\widehat{\alpha} \circ \widehat{\rho}_i= \widehat{\rho}_{\alpha(i)} \circ \widehat{\alpha}$ for every
$\alpha\in\mathfrak{S}_{2m}$.

\end{enumerate}
\end{prop}

\subsection{Generalized Kashaev algebra: triangulation-independent}

\begin{theorem}[Guo-Liu \cite{GL}]\label{thm:main}
For two arbitrary complex numbers $a,b$, there is a family of algebra isomorphisms
$$\Psi_{\tau\tau'}^q(a,b):
\widehat{\mathcal{K}}^q_{\tau'} \rightarrow
\widehat{\mathcal{K}}^q_{\tau}$$ defined as $\tau$,
$\tau' \in \triangle(S)$ ranges over all pairs of decorated ideal
triangulations, such that:
\begin{enumerate}
\item $\Psi_{\tau\tau''}^q(a,b) = \Psi_{\tau\tau'}^q(a,b) \circ
\Psi_{\tau'\tau''}^q(a,b)$ for every $\tau$, $\tau'$,
$\tau''\in \triangle(S)$;

\item $\Psi_{\tau\tau'}^q(a,b)$ is the isomorphism of
Definition \ref{def:iso} when $\tau'$ is obtained
from $\tau$ by a reindexing or a mark rotation or a diagonal exchange.

\item $\Psi_{\tau\tau'}^q(a,b)$ depends only on $\tau$ and $\tau'$.
\end{enumerate}
\end{theorem}

The \emph{generalized Kashaev algebra}\index{Kashaev algebra!generalized} $\widehat{\mathcal{K}}^q_S(a,b)$ associated to a surface $S$ is defined as the algebra
$$
\widehat{\mathcal{K}}^q_S(a,b)= \bigg(
\bigsqcup_{\tau \in \triangle(S)}
\widehat{\mathcal{K}}^q_\tau(a,b)\bigg)/\sim
$$
where the relation $\sim$ is defined by the property that, for
$X\in \widehat{\mathcal{K}}^q_\tau(a,b)$ and $X'\in
\widehat{\mathcal{K}}^q_{\tau'}(a,b)$,
$$
X \sim X' \Leftrightarrow X=\Psi^q_{\tau,\tau'}(a,b)(X').
$$

\section{Kashaev coordinates and shear coordinates}

To understand the quantization using shear coordinates and the quantization using Kashaev coordinates, we first need to understand the relationship between these two coordinates.

\subsection{Decorated ideal triangulations}
Given a decorated ideal triangulation $\tau,$ by forgetting the mark at each corner, we obtain an ideal triangulation $\lambda.$ We call $\lambda$ the \emph{underlying ideal triangulation}\index{ideal triangulation!underlying} of $\tau$. let $\lambda_1,\lambda_2,...,\lambda_{3m}$ be the components of ideal triangulation $\lambda$. Denote by $\tau_1,..,\tau_{2m}$ the ideal triangles.

For the ideal triangulation $\lambda$, we may consider its dual graph. Each ideal triangle $\tau_\mu$ corresponds to a vertex $\tau_\mu^*$ of the dual graph. Denote by $\lambda_1^*,\lambda_2^*,...,\lambda_{3m}^*$ the dual edges. If an edge $\lambda_i$ bounds one side of the ideal triangles $\tau_\mu$ and one side of $\tau_\nu$, then the dual edge $\lambda_i^*$ connects the vertexes $\tau_\mu^*$ and $\tau_\nu^*$.

In a decorated ideal triangulation $\tau$, each ideal triangle $\tau_\mu$ (embedded or not) has three sides which correspond to the three half-edges incident to the vertex $\tau_\mu^*$ of the dual graph. The three sides are numerated by $0,1,2$ in the counterclockwise order such that the $0-$side is opposite to the marked corner.

\subsection{Space of Kashaev coordinates}
Let's recall that a Kashaev coordinate associated to a decorated ideal triangulation $\tau$ is a vector $(\ln y_1,\ln z_1,...,\ln y_{2m},\ln z_{2m})\in \mathbb{R}^{4m}$, where $\ln y_\mu$ and $\ln z_\mu$ are associated to the ideal triangle $\tau_\mu$. Denote by $\mathcal K_\tau$ the space of Kahsaev coordinates associated to $\tau$. We see that $\mathcal K_\tau=\mathbb{R}^{4m}.$

Given a vector $(\ln y_1,\ln z_1,...,\ln y_{2m},\ln z_{2m})\in \mathcal K_\tau$, we associate a number to each side of each ideal triangle as follows. For the ideal triangle $\tau_\mu$, we associate
\begin{itemize}
\item $\ln h_\mu^0:=\ln y_\mu-\ln z_\mu$ to the 0-side;

\item $\ln h_\mu^1:=\ln z_\mu$ to the 1-side;

\item $\ln h_\mu^2:=-\ln y_\mu$ to the 2-side.
\end{itemize}

Therefore $\ln h_\mu^0+\ln h_\mu^1+\ln h_\mu^2=0$. We can identify the space  $\mathcal K_\tau=\mathbb{R}^{4m}$ with a subspace of $\mathbb{R}^{6m}=\{(...,\ln h_\mu^0,\ln h_\mu^1,\ln h_\mu^2,...)\}$ satisfying $\ln h_\mu^0+\ln h_\mu^1+\ln h_\mu^2=0$ for each ideal triangle $\tau_\mu$.

\subsection{Exact sequence}
The enhanced Teichm\"uller space parametrized by shear coordinates is $\widetilde{\mathcal T}_\lambda=\mathbb{R}^{3m}=\{(\ln x_1,\ln x_2,...,\ln x_{3m})\}$, where $\ln x_i$ is the shear coordinate at edge $\lambda_i$. We define a map $f_1: \widetilde{\mathcal T}_\lambda \to \mathbb{R}$ by sending $(\ln x_1,\ln x_2,...,\ln x_{3m})$ to the sum of entries $\sum_{i=1}^{3m}\ln x_i.$

Suppose $\lambda$ is the underlying ideal triangulation of the decorated ideal triangulation $\tau$. We define a map $f_2: \mathcal K_\tau \to \widetilde{\mathcal T}_\lambda$ as a linear function by setting $$\ln x_i=\ln h_\mu^s+\ln h_\nu^t$$ whenever $\lambda_i$ bounds the $s-$side of $\tau_\mu$ and the $t-$side of $\tau_\nu$ ($\mu$ may equal $\nu$), where $s,t\in\{0,1,2\}$.

Another map $f_3: H_1(S,\mathbb{R}) \to \mathcal K_\tau$ is defined as follows. A homology class in $H_1(S,\mathbb{R})$ is represented by a linear combination of oriented dual edges: $\sum_{i=1}^{3m}c_i\lambda_i^*$. If the orientation of $\lambda_i^*$ is from the $s-$side of $\tau_\mu$ to the $t-$side of $\tau_\nu$, by setting $\ln h_\mu^s=-c_i$ and $\ln h_\nu^t=c_i$, we obtain a vector $(...,\ln h_\mu^0,\ln h_\mu^1,\ln h_\mu^2,...)\in \mathbb{R}^{6m}.$ The boundary map of chain complexes sends $\sum_{i=1}^{3m}c_i\lambda_i^*$ to a linear combination of vertexes. In this combination, the term involving the vertex $\tau_\mu^*$ is $(c_i\epsilon_i +c_j\epsilon_j+c_k\epsilon_k)\tau_\mu^*$ where $\lambda_i,\lambda_j,\lambda_k$ (two of them may coincide) bound three sides of $\tau_\mu$ and $\epsilon_t=-1$ if $\lambda_t^*$ starts at the side of $\tau_\mu$ bounded by $\lambda_t$ while $\epsilon_t=1$ if $\lambda_t^*$ ends at the side of $\tau_\mu$ bounded by $\lambda_t$. Therefore $$(c_i\epsilon_i +c_j\epsilon_j+c_k\epsilon_k)\tau_\mu^*=(\ln h_\mu^0+\ln h_\mu^1+\ln h_\mu^2)\tau_\mu^*.$$ Since the chain $\sum_{i=1}^{3m}c_i\lambda_i^*$ is a cycle, we must have $\ln h_\mu^0+\ln h_\mu^1+\ln h_\mu^2=0$. Therefore this vector $(...,\ln h_\mu^0,\ln h_\mu^1,\ln h_\mu^2,...)$ is in the subspace $\mathcal K_\tau$.

Combining the three maps, we obtain

\begin{theorem}[Guo-Liu \cite{GL}]\label{thm:exact} The following sequence is exact:
$$0\rightarrow H_1(S, \mathbb R)\xrightarrow{f_3} \mathcal K_\tau \xrightarrow{f_2} \widetilde{\mathcal T}_\lambda \xrightarrow{f_1} \mathbb R \rightarrow 0.$$
\end{theorem}

From the theorem above, we see that $\mathcal K_\tau$ is a fiber bundle over the space $\mathrm{Ker}(f_1)$ whose fiber is an affine space modeled on $H_1(S, \mathbb R).$ To be precise, given a vector $s\in \mathrm{Ker}(f_1),$ let $v\in f_2^{-1}(s).$ Then $f_2^{-1}(s)=v+H_1(S, \mathbb R).$

\subsection{Relation to bivecotrs}

Consider the linear isomorphism
\begin{align}\label{fml:3-2}
M:  \mathcal K_\tau &\longrightarrow  \mathcal K_\tau\\
(\ln y_1,\ln z_1,...,\ln y_{2m},\ln z_{2m})&\longmapsto  (...,\ln h_\mu^0,\ln h_\mu^1,\ln h_\mu^2,...) .\notag
\end{align}

\begin{prop}[Guo-Liu \cite{GL}]\label{thm:bivector}~

\noindent If $(\ln x_1,\ln x_2,...,\ln x_{3m})=f_2\circ M (\ln y_1,\ln z_1,...,\ln y_{2m},\ln z_{2m}),$ then
$$\sum_{i,j=1}^{3m}\sigma^\lambda_{ij}\frac{\partial}{\partial \ln x_i} \wedge \frac{\partial}{\partial \ln x_j}=
(f_2)_*\circ M_* (\sum_{\mu=1}^{2m}\frac{\partial}{\partial \ln y_\mu} \wedge \frac{\partial}{\partial \ln z_\mu}),$$
where $\sigma^\lambda_{ij}=a^\lambda_{ij}-a^\lambda_{ji}$ and $a^\lambda_{ij}$ is the number of corners of the ideal triangulation $\lambda$ which is delimited in the left by $\lambda_i$ and on the right by $\lambda_j$.
\end{prop}

The left hand side of the equality is the Weil-Petersson Poisson structure\index{Weil-Petersson Poisson structure} on the enhanced Teichm\"uller space \cite{FC}.

\subsection{Compatibility of coordinate changes}

\begin{prop}[Guo-Liu \cite{GL}]\label{compatible} Suppose that the decorated ideal triangulations $\tau$ and $\tau'$ have underlying ideal triangulations $\lambda$ and $\lambda'$ respectively. The following diagram is commutative:
$$
\begin{CD}
\widetilde{\mathcal{T}}_\lambda @<f_2<< \mathcal{K}_\tau\\
@VVV @VVV \\
\widetilde{\mathcal{T}}_{\lambda'} @<f_2<<  \mathcal{K}_{\tau'}
\end{CD}
$$
where the two vertical maps are corresponding coordinate changes. The coordinate changes of Kashaev coordinates are given in Definition \ref{def:coor-change}. The coordinate changes of shear coordinates are given in Proposition \ref{prop:DiagExch}.
\end{prop}

\section{Relationship between quantum Teichm\"uller space and Kashaev algebra}

In this section, we establish a natural relationship between the quantum Teichm\"uller space $\widehat{\mathcal{T}}^q_S$ and the generalized Kashaev algebra $\widehat{\mathcal{K}}^q_S(a,b)$.

\subsection{Homomorphism}

For a ideal triangle $\tau_\mu$, we associate three elements in $\mathcal{K}^q_{\tau}$ to the three sides of $\tau_\mu$ as follows:
\begin{itemize}
\item $H^0_\mu:=Y_\mu Z^{-1}_\mu$ to the $0-$side;

\item $H^1_\mu:=Z_\mu$ to the $1-$side;

\item $H^2_\mu:=Y^{-1}_\mu$ to the $2-$side.
\end{itemize}

\begin{lemma}\label{thm:h} For any $s,t\in\{0,1,2\}$ and $\mu\in {1,2,...,3m},$
$$H^s_\mu H^t_\mu=q^{2\sigma_{st}}H^t_\mu H^s_\mu,$$
where $\sigma_{st}+\sigma_{ts}=0$ and $\sigma_{10}=\sigma_{02}=\sigma_{21}=1.$
\end{lemma}

Suppose $\lambda$ is the underlying ideal triangulation of $\tau$, the Chekhov-Fock algebra $\mathcal T_\lambda^q$ is the algebra over $\mathbb C$ defined by generators $X_1^{\pm1}$, $X_2^{\pm1}$, \dots, $X_n^{\pm1}$ associated to the components of $\lambda$ and by relations $X_iX_j=q^{2\sigma^\lambda_{ij}}X_jX_i$.

We define a map $F_\tau$ from the set of the generators of $\mathcal T_\lambda^q$ to $\mathcal K^q_{\tau}$. Suppose that the edge $\lambda_i$ bounds the $s-$side of $\tau_\mu$ and the $t-$side of $\tau_\nu$. We define

\begin{align}\label{def:X}
F_\tau(X_i)=q^{\delta_{\mu\nu}\sigma_{ts}}H^s_\mu H^t_\nu \in \mathcal{K}^q_\tau,
\end{align}
where $\sigma_{ts}$ is defined in Lemma \ref{thm:h} and $\delta_{\mu\nu}$ is the Kronecker delta, i.e., $\delta_{\mu\mu}=1$ and $\delta_{\mu\nu}=0$ if $\mu\neq \nu$. When $\mu=\nu,$ $X_i$ is well-defined, since $$q^{\sigma_{ts}}H^s_\mu H^t_\mu=q^{\sigma_{st}}H^t_\mu H^s_\mu$$ due to Lemma \ref{thm:h}.

This definition is natural since when $q=1$ we get the relationship between Kashaev coordinates and shear coordinates which is given by the map $M$ and $f_2$. In fact when $q=1$ the generators $Y_\mu, Z_\mu$ are commutative. They reduce to the geometric quantities $y_\mu,z_\mu$ associate to $\tau_\mu.$ $H^s_\mu$ and $X_i$ are reduced to $h^s_\mu$ and $x_i$.

\begin{lemma} $F_\tau(X_i)F_\tau(X_j)=q^{2\sigma^\lambda_{ij}}F_\tau(X_j)F_\tau(X_i)$ for any generators $X_i$ and $X_j$.
\end{lemma}

It turns out that $F_\tau$ can be extended to the whole algebra $\mathcal T_\lambda^q$ as an algebra homomorphism
from $\mathcal T_\lambda^q$ to $\mathcal K^q_{\tau}$.

\subsection{Compatibility}

Recall that $\widehat{\mathcal K}^q_{\tau}$ is the fraction division algebra of $\mathcal K^q_{\tau}$. The algebraic isomorphism between $\widehat{\mathcal K}^q_{\tau}$ and $\widehat{\mathcal K}^q_{\tau'}$ is defined in Definition \ref{def:iso}.

\begin{lemma}\label{thm:a}
Suppose that a decorated ideal triangulation $\tau'$ is obtained from $\tau$ by a mark rotation $\rho_\mu$ for some $\mu\in\{1,2,...,2m\}$. Let $\lambda$ be the common underlying ideal triangulation of $\tau$ and $\tau'$. The following diagram is commutative if and only if $a=q^{-2}$.
$$
\begin{CD}
\widehat{\mathcal T}_\lambda^q @>F_\tau>> \widehat{\mathcal K}^q_{\tau}\\
@A\mathrm{Id}AA @AA\widehat{\rho}_\mu A \\
\widehat{\mathcal T}_\lambda^q @>F_{\tau'}>>  \widehat{\mathcal K}^q_{\tau'}
\end{CD}
$$

\end{lemma}

\begin{lemma}\label{thm:b}
Suppose that a decorated ideal triangulation $\tau'$ is obtained from $\tau$ by a diagonal exchange $\varphi_{\mu\nu}.$ Let $\lambda$ and $\lambda'$ be the underlying ideal triangulation of $\tau$ and $\tau'$ respectively. Then $\lambda'$ is obtained from $\lambda$ by a diagonal exchange with respect to the edge $\lambda_i$ which is the common edge of $\tau_\mu$ and $\tau_\nu$.
The following diagram is commutative if and only if $b=q^{3}$.
$$
\begin{CD}
\widehat{\mathcal T}_\lambda^q @>F_\tau>> \widehat{\mathcal K}^q_{\tau}\\
@A\widehat{\Delta}_i AA @AA\widehat{\varphi}_{\mu\nu} A \\
\widehat{\mathcal T}_{\lambda'}^q @>F_{\tau'}>>  \widehat{\mathcal K}^q_{\tau'}
\end{CD}
$$
\end{lemma}

\begin{theorem}[Guo-Liu \cite{GL}]
Suppose the decorated ideal triangulations $\tau$ and $\tau'$ have underlying ideal triangulations
$\lambda$ and $\lambda'$ respectively.
The following diagram is commutative if and only if $a=q^{-2}, b=q^{3}$.
$$
\begin{CD}
\widehat{\mathcal T}_\lambda^q @>F_\tau>> \widehat{\mathcal K}^q_{\tau}\\
@A\Phi^q_{\lambda,\lambda'} AA @AA\Psi^q_{\tau,\tau'}(a,b) A \\
\widehat{\mathcal T}_{\lambda'}^q @>F_{\tau'}>>  \widehat{\mathcal K}^q_{\tau'}
\end{CD}
$$
\end{theorem}

Recall that the quantum Teichm\"uller space of $S$ is defined as the algebra
$$
\widehat {\mathcal{T}}^q_S= \bigg(
\bigsqcup_{\lambda\in\Lambda(S)}
\widehat{\mathcal{T}}^q_{\lambda}\bigg)/\sim
$$
where the relation $\sim$ is defined by the property that, for
$X\in \widehat{\mathcal{T}}^q_{\lambda}$ and $X'\in
\widehat{\mathcal{T}}^q_{\lambda'}$,
$$
X \sim X' \Leftrightarrow X=\Phi^q_{\lambda,\lambda'}(X').
$$

The generalized Kashaev algebra $\widehat{\mathcal{K}}^q_S(a,b)$ associated to a surface $S$ is defined as the algebra
$$
\widehat{\mathcal{K}}^q_S(a,b)= \bigg(
\bigsqcup_{\tau \in \triangle(S)}
\widehat{\mathcal{K}}^q_\tau(a,b)\bigg)/\sim
$$
where the relation $\sim$ is defined by the property that, for
$X\in \widehat{\mathcal{K}}^q_\tau(a,b)$ and $X'\in
\widehat{\mathcal{K}}^q_{\tau'}(a,b)$,
$$
X \sim X' \Leftrightarrow X=\Psi^q_{\tau,\tau'}(a,b)(X').
$$

\begin{corollary}\label{cor:homo} The homomorphism $F_\tau$ induces a homomorphism $$\widehat {\mathcal{T}}^q_S \to \widehat{\mathcal{K}}^q_S(a,b)$$ if and only if $a=q^{-2}, b=q^3$.
\end{corollary}

\subsection{Quotient algebra}

Furthermore, consider the element
$$H=q^{-\sum_{i<j}\sigma^\lambda_{ij}}X_1X_2...X_{3m}\in \mathcal T_\lambda^q.$$ It is proved in \cite{Liu1}(Proposition 14) that $H$ is independent of the ideal triangulation $\lambda$. Therefore $H$ is a well-defined element of the
quantum Teichm\"uller space $\widehat{\mathcal{T}}^q_S$.

\begin{theorem}[Guo-Liu \cite{GL}]\label{thm:homo} The homomorphism $F_\tau$ induces a homomorphism $$\widehat {\mathcal{T}}^q_S/(q^{-2m}H-1)\to \widehat{\mathcal{K}}^q_S(q^{-2},q^3)$$ where $(q^{-2m}H-1)$ is the ideal generated by $q^{-2m}H-1$.
\end{theorem}

\end{document}